\documentclass[a4paper,11pt]{amsart}

\usepackage{fancyhdr}
\usepackage{amsmath,inputenc,euscript,amssymb}
\usepackage{fullpage}
\usepackage{bbm}
\usepackage{color}
\usepackage[colorlinks=true, linkcolor=magenta, citecolor=cyan, urlcolor=blue]{hyperref}
\usepackage{datetime}

\font\bb=msbm10 at 11pt
\font\bbi=msbm10 at 8pt

\def\R{\hbox{\bb R}}
\def\C{\hbox{\bb C}}

\def\Ri{\hbox{\bbi R}}

\newtheorem{theo}{Theorem}[section]
\newtheorem{pro}[theo]{Proposition}
\newtheorem{lem}[theo]{Lemma}

\newtheorem{cor}[theo]{Corollary}

\newtheorem{rem}[theo]{Remark}

\title{Revisiting the Blow-up Criterion and the Maximal Existence Time  for Solutions 
of the Parabolic-Elliptic Keller-Segel System in 2D-Euclidean Space}



\author{Patrick Maheux \&  Vittoria Pierfelice}
\address{Institut Denis Poisson (IDP)\\
Universit\'e d'Orl\'eans, Universit\'e  de Tours, CNRS Ð UMR 7013
\\ B\^atiment de Math\'ematiques, B.P. 6759\\
45067 Orl\'eans Cedex 2, France}

\email{Patrick.Maheux@univ-orleans.fr
\& Vittoria.Pierfelice@univ-orleans.fr}

\begin{document}

\maketitle

\begin{center}
\today  
\end{center}

\begin{abstract}
In this paper, we revisit the blow-up criteria for the simplest parabolic-elliptic Patlak–Keller–Segel (PKS) system in the 2D Euclidean space, including a consumption term.
In the supercritical mass case $M > 8\pi$, and under an additional global assumption on the second moment (or variance) of the initial data, we establish blow-up results for a broader class of initial conditions than those traditionally considered.
We also derive improved upper bounds for the maximal existence time of (PKS) solutions on the plane. These time estimates are obtained through a sharp analysis of a one-parameter differential inequality governing the evolution of the second moment of the (PKS) system.
As a consequence, for any given non-negative (non-zero) initial datum $n_1$ with finite second moment, we construct blow-up solutions of the (PKS) system.
\end{abstract}

\tableofcontents

\section{Introduction} \label{intro}

The aim of this paper is to demonstrate that enlarge  the currently known  upper bounds  on the second moment
of the  initial data,  leading to blow-up of the solutions for the Patlak-Keller-Segel  system  on the plane  (considered  with a consumption term) 
for the  supercritical mass  case $M>8\pi$.

The goal of the present study for  this chemotaxis model  is twofold. 
First, we search for a differential inequality satisfied by the variance of the solution of the (PKS) system \eqref{ks} in terms of a function 
$g_{\alpha}$ intrinsically related to this equation by using a convexity argument instead of approximations of $g_{\alpha}$.
The second aspect consists in solving the differential inequality \eqref{diffit2} in an optimal way to get a precise criterion on the critical variance 
$V(0)$ at the initial time $t=0$ for blow-up and also an upper bound on the blow-up time $T^*$ of the solutions of (PKS) system. 

More precisely, the one-parameter differential  inequality obtained for the variance  $V(t)$ (or moment $I(t)$)  function is  solved  in a  general form
for itself, i.e. independently of the (PKS) system from which it comes from. 
More specifically, there are  two aspects: 
\vskip0.1cm
 (i) finding the largest  upper  bound  on  the initial  data  $V(0)$, which ensures the blow-up of the solution  $(V(t), 0<t<T^*)$,
  i.e.  $T^*<+\infty$, 
for   this type of differential inequality.
\vskip0.1cm
(ii) Obtaining  an optimal upper  bound  on  the maximal existence  time  $T^*$ for  the class of solutions $(V(t), 0<t<T^*)$
satisfying  this differential  inequality.
\vskip0.3cm

We then apply these results  to  the evolution in time of the  variance $V(t)$ of the solution of the (PKS) system.
This  yields a blow-up criterion based on the initial variance $V(0)$, and  also an upper-bound on the maximal existence  time  for this solution.
We note that the notion of variance deserves further emphasis because 
it encompasses both the second moment and the center of mass  in a single expression.
Therefore, in this work, we  express our results  in terms of  variance.
\vskip0.3cm

The whole approach  unifies   and improves several proofs of blow-up criteria 
and bounds of  maximal existence  time  of the current  literature
for  the (PKS) system considered in this paper.

\vskip0.3cm
{\em The content of this paper is the following}.
\vskip0.3cm

 In Section \ref{intro}, 
we introduce  the well-known Patlak-Keller-Segel  (simplest) system  with consumption term on $\R^2$  
considered in this paper. In particular, we discuss the role of  the Bessel kernel and  the associated function  $g_{\alpha}$
 important for the study of this system.
(Section \ref{PKSintro}).
  
We also recall  known results 
about  blow-up  criteria on the second moment of the initial data $n_0$ of the (PKS) 
solutions obtained in \cite{CC1} and \cite{koz-sug}.
We also present  precise upper-bounds  on the blow-up time $T^*$,
i.e.  bounds on the maximal existence   time  $T^*$
of  the solutions of  the (PKS) system, presented in the aforementioned papers. (Section \ref{known}).
\vskip0.3cm
  
In Section \ref{mainresult},
 we  first comment on  the previous results recalled in Section \ref{known}, which
motivate our work. We also briefly compare these results with some  of those presented  here.   
(Section \ref{aim}).
  
We state the main result (Theorem \ref{bupp})  by  introducing a new (global)
 blow-up  criterion on the   variance (or second moment)  of the initial data.
Simultaneously,   we obtain precise  bounds 
on the maximal time of existence of the solutions of the (PKS) system.
The main  feature    
of our result is that  the blow-up criterion  is expressed  in terms of  the function  $g_{\alpha}$,
described  by  \eqref{galphafn} below,  directly related to the Bessel function $B_{\alpha}$, defined  in \eqref{besselkern}, but 
without using  any estimate of this function $g_{\alpha}$. 
This provides a result of  an intrinsic nature  for the blow-up criterion.
 The key argument  consists in   applying Jensen's inequality appropriately  by  noting the convexity of  the function
$\rho\rightarrow g_{\alpha}(\sqrt{\rho})$. This leads to an elegant  and natural proof of  the blow-up criterion  for the (PKS) system.
(Section \ref{result}).

Independently, we develop a complete study of a certain  
family of one-parameter  differential inequalities (Proposition \ref{diffineq}).
In particular, this study  provides a blow-up criterion  on the initial data  and a  sharp bound  of 
 the maximal time of existence  of  the solutions. These results are  applied specifically to the (PKS) system treated  in this paper.  
(Section \ref{diff-argument}).
 \vskip0.3cm
       
 In Section \ref{lescoro}, we derive several corollaries of Theorem \ref{bupp}.
 In Section \ref{deuxcoro},  we discuss the behavior of the blow-up  time  of the 
 solutions when considering multiples $\lambda n_0$ of  fixed initial data $n_0$.
 The first result,
Corollary \ref{timedecay},    deals with the time decay in terms of $\lambda\geq 1$   when
the  initial data $n_0$ satisfies the blow-up criterion.
More interestingly,  the second result, Corollary \ref{blowcreation}, allows us  to construct examples of blow-up solutions
from any initial data  with  finite 2-moment by taking large multiples  $\lambda$ of it.
This result seems quite new to the authors' knowledge.

In Section \ref{lastcoro}, 
we derive several Corollaries  \ref{firstcoro}, \ref{secondcoro} and  \ref{explicitbup} 
using  estimates of  the functions $g_{\alpha}$ and  $g^{-1}_{\alpha}$ described in Proposition \ref{propgun} 
of  Section \ref{gestim} to make the blow-up condition  more explicit.

In  Section \ref{gestim}, additional  estimates on $g_{\alpha}$ allows us to  precisely compare 
the growth of the previous  bounds with the new bound imposed on the initial variance
(Proposition \ref{progamcc}). 
This comparison  leads  to a larger  blow-up class of initial data.

\subsection{The Patlak-Keller-Segel (PKS) model for chemotaxis}\label{PKSintro}
In this paper, we consider the following simple Patlak-Keller-Segel system  
\eqref{pks1}-\eqref{pks2}, (PKS) for short, given below and 
modellling a phenomenon of  chemotaxis on  the whole plane:
\begin{equation}\label{pks1}
\frac{\partial}{\partial t}n(t,x)=\Delta n(t,x)-\nabla.(n\nabla c)(t,x), \quad t>0, \; x\in \R^2,
\end{equation}
\begin{equation}\label{pks2}
(-\Delta)c+\alpha c=n, \quad t>0, \; x\in \R^2,
\end{equation}
with 
$n(0,x)=n_0(x)$, $0\leq n_0\in L^1(\R^2)$ as initial data, and $\alpha\geq 0$ a fixed parameter.
See \cite{KellerSegel},\cite{Patlak}.
\vskip0.3cm

Here, the letter $n$ denotes the density of a cell population and $c$ is the concentration of a chemical signal attracting the cells. 
The constant $\alpha \geq 0$ is called  the consumption term or chemical degradation rate. 
See \cite{CC1}, for instance.

More precisely, the (PKS)  system describes the collective motion of cells which
are attracted by a chemical substance and are able to emit it. 
See \cite{KellerSegel, Patlak} for the origin for such a  system,  
and the review papers \cite{Hort1, Hort2}.
We can also consult \cite{Mur} for a general introduction to chemotaxis
and  the recent book by P.Biler \cite{Bil} for some account of more recent results on  (PKS)  systems. 
See also the  references therein.
\vskip0.3cm

Many  results  of local or global   existence in time of the solution
$n_t(x)=n(t,x)$, $x\in \R^2$, of the (PKS) system have been proved with the initial conditions:
$n(0,x)=n_0(x)\geq 0$ and  $n_0\in L^1(\R^2)$ where  additional conditions on $n_0$ are imposed. For instance, 
finite second moment  and  finite entropy  assumptions on the initial data  are requested
(or, more generally in a weighted $L^1(wdx)$ for some weights $w$).
Recently in  \cite{DongyiWei} (for the case $\alpha=0$),
only the assumptions of integrability and nonnegativity of $n_0$  are imposed for the existence of mild solutions.
No finite second moment condition or  finite entropy  on $n_0$ is required.
In particular, the author of  \cite{DongyiWei}  proves  that a global in time solution exists  if  the mass $M\leq 8\pi$
and  the (local) solution blows up if  $M >8\pi$. 
\vskip0.3cm

The (PKS) system can be  interpreted mathematically in different ways, and  solved in different functional spaces of solutions.
Consequently, the solutions are called by  different names : strong, classical, weak,  mild solutions, etc... 
In this paper,  we   will   assume more or less implicitly  that our computations  concerning the blow-up criterion  can be justified 
as in \cite{CC1}, or in  \cite{koz-sug}. 
But in the setting of  \cite{CC1} for instance,  we   will   not use explicitly of  the   entropy assumption  
on the initial data  $n_0$ in our computations (only finite second moment is required).
For results of blow-up criteria  using finite entropy of the initial data in higher dimension, we can consult  \cite[Prop.4.2]{CC3}  for instance.
Note that the blow-up criteria  in  \cite{CC1} and   \cite{koz-sug}, as  in this paper, are of  global type, 
i.e.  the whole space is involved in the criterion.  More recently and more natural, local criteria 
on the initial data have been used for blow-ups, see \cite{BCKZ,BKZ}.
\vskip0.3cm

There is an abundant  literature on the study  (PKS) system, its variants or generalizations.
We apologize for mentioning not a complete bibliography on the subject for this reason.
 We just mention a few of them in this paper concerning simple systems  on $\R^2$ with $\alpha=0$ or $\alpha \geq 0$, which have some direct links to our study
 (but this list is by no means exhaustive):
  \cite{BCKZ},
 \cite{CC1},
\cite{Ka-Su}, 
 \cite{koz-sug},
 \cite{Nagai-0},
 \cite{Nagai4},
 \cite{Nagai4-0},
\cite{Nagai4-1},   
\cite{Mizo},  
\cite{Sen1},
\cite{DP}.
See also
\cite{Jager},
\cite{HV1},
for related problems or former results.
General references on these issues  are:
\cite{Hort1},\cite{Hort2},\cite{Perth}.
Recently, a result of blow-up has been obtained in \cite{Ma-Pier} on the hyperbolic space of dimension two.

In  this paper, we assume  that the initial condition $n_0\in L^1(\R^2)$ is nonnegative and non zero, and has a finite second moment.
We consider a setting where the local existence of the solution is known   (plus additional conditions if needed).
We   will   not discuss existence issues  in this paper,
but we focus on blow-up criteria and bounds on the maximal existence time of the solutions.
In the whole paper, $T^*$   will   denote  the maximal time of existence of the (mild, weak, strong,...) 
solution $(n_t)$, $0<t<T^*\leq +\infty$, of the (PKS) system.
\vskip0.3cm

Now, we introduce notions  that  are involved in our considerations of the solutions of (PKS) systems for $\alpha >0$ as done in \cite{CC1}.
The case $\alpha =0$ is well studied, see \cite{BDP,DP,DongyiWei, Bil} for instance.
\vskip0.3cm

In the whole paper, we   will   only consider the solution $c$ given by  
$$
c(t,x)=(-\Delta+\alpha)^{-1} n(t,x), \; x\in \R^2,\, t>0,
$$
 obtained by {\it shifting} the spectral gap of the (non-negative) Laplacian $-\Delta$ by $\alpha >0$.
Thus, we first  focus our analysis on the resolvent operator $(-\Delta+\alpha)^{-1}$ of the Laplacian $-\Delta$.
More precisely, we consider the solution $c_t$ given by the  convolution with  the  kernel $B_{\alpha}$ and $n_t$, i.e.  
$
c_t=
(-\Delta+\alpha)^{-1}n_t=
B_{\alpha}\ast n_t$,
where 
\begin{equation}\label{besselkern}
B_{\alpha}(z)=\frac{1}{4\pi}\int_0^{+\infty} \frac{1}{t}e^{\frac{-\vert z\vert^2}{4t}}
e^{-\alpha t}\
\, dt
={\mathcal L}(K_{.}(z))(\alpha).
\end{equation}
(Here, $\ast$ denotes the convolution of functions on $\R^2$).  The kernel $B_{\alpha}$ is the Bessel kernel obtained as
 the Laplace transform ${\mathcal L}(K_{.}(z))(\alpha)$   in time
of  the two-dimensional  heat kernel $K_t$ on $\R^2$,  defined  by
$$
K_t(z)= \frac{1}{4\pi t} e^{\frac{-\vert z\vert^2}{4t}}, \qquad t>0, \;  z\in \R^2.
$$

Due to  the facts that the solution $n_t\geq 0$ and $B_{\alpha}\geq 0$, we have $c_t\geq 0$.
\vskip0.3cm

From the discussion just above, the (PKS) system  \eqref{pks1}-\eqref{pks2} above  reduces to  the following   nonlinear parabolic equation,
\begin{equation}\label{ks}
\frac{\partial}{\partial t}n(t,x)=\Delta n(t,x)-\nabla.(n\nabla(B_{\alpha}\ast n ))(t,x), \; t>0, \, x\in \R^2.
\end{equation}
We   will   call this equation the (PKS) equation.

The gradient $\nabla B_{\alpha}$ of the kernel $B_{\alpha}$ plays an important role  in the nonlinear term of the equation \eqref{ks} 
for a possible blow-up.
This gradient  is related to the next function $g_{\alpha}$ ($\alpha>0$), 
crucial for our analysis and  defined  by
\begin{equation}\label{galphafn}
g_{\alpha}(r)=\int_0^{+\infty}e^{-\frac{\alpha r^2}{4s}}  e^{-s}\, ds, \quad r\geq 0,
\end{equation}
via   the next formula
\begin{equation}\label{linkbesselgalpha}
\nabla B_{\alpha}(z)
=
-\frac{z}{2\pi \vert z\vert^2}
\int_0^{+\infty}
e^{-s} e^{-\alpha \frac{\vert z\vert^2}{4s}}\,ds
=
-
\frac{z}{2\pi \vert z\vert^2} \, g_{\alpha}(\vert z\vert),\; z\in \R^2.
\end{equation}
 (Our notation of $g_{\alpha}$  is slightly different from  the usual one.
Here, the   expression  of $g_{\alpha}(\vert x-y\vert)$ corresponds to 
 $g_{\alpha}(x-y)$ in  \cite[Section 7]{CC1}).

In the literature, various estimates of the  function $g_{\alpha}$  lead to explicit  blow-up criteria, 
i.e. explicit upper bounds on the second moment $I(0)$  of the initial data  $n_0$ with 
\begin{equation} 
I(0):=
\int_{\Ri^2} \vert x\vert^2\, n_0(x)dx.
\end{equation}
Bounds of the maximal time of existence $T^*$ of  the solutions of the (PKS) equation \eqref{ks} are also obtained.
 See for instance  \cite{CC1}, \cite{koz-sug} and our study below.
\vskip0.3cm

\subsection{Known results on the blow-up with consumption term}\label{known}
From now on, we assume that the initial data  $n_0=n(0,.)$ is non-negative and non-zero with  finite  second moment,  i.e.
\begin{equation}\label{secondmoment}
I(0):=\int_{\Ri^2} \vert x\vert^2n_0(x)\, dx<+\infty.
\end{equation}
To study the blow-up  of the  solution of the (PKS) equation \eqref{ks} with initial condition $n_0$, we  classically consider  the evolution in time of the 
second moment of the solution defined by
\begin{equation}\label{ioft}
t\in (0,T^*) \mapsto I(t)=
\int_{\Ri^2} \vert x\vert^2n(t,x)\,dx,
\end{equation}
where $T^*$ is the maximal  existence time of the solution $(n_t)$. 
Recall that the total mass $M$ is constant during the time evolution, i.e.
$$
M:=
\int_{\Ri^2}  n_0(x)\,dx=\int_{\Ri^2}  n(t,x)\,dx  , \quad  0<t<T^*.
$$

In \cite[Theorem 1.2]{CC1}, it has been proved that if the following two conditions $M>8\pi$, and
\begin{equation}\label{ccond}
I(0)=\int_{\Ri^2} \vert x\vert^2n_0(x)\, dx< \mu^*,
\end{equation}
are satisfied with 
\begin{equation}\label{lamc}
\mu^*=\mu^*(\alpha,M)= \frac{1}{\alpha} h_1(M/8\pi)
=
\frac{1}{4\alpha \, \mathcal{C}^2 M}
(M-8\pi)^2,
\end{equation}
where $h_1(s)=\frac{2\pi}{ { \mathcal{C}}^2}.\frac{(s-1)^2}{s}$, $s>1$,  and for some constant $  \mathcal{C}\geq 1$,
then the solution $(n_t)$ blows up in finite time $T^*$
such that 
\begin{equation}\label{tcc}
T^*\leq T^*_{cc}(\alpha,M, I(0)):=
\frac{ I(0)}
{4M(\frac{M}{8\pi}-1) -\frac{{\mathcal C}}{\pi}  {\alpha}^{1/2} M^{3/2} I(0)^{1/2}}.
\end{equation}
For details on the constant $\mathcal C\geq 1$ and the proof, see \cite[Section 7]{CC1}.
See also \cite{Jun-Lein} for the use of these bounds $\mu^{*}$ and $T^*_{cc}$ to the study of discret (PKS)  models.

Around the same period of \cite{CC1}, the authors of \cite[Th2.p.356]{koz-sug}  proved that if,
instead of \eqref{ccond}, the condition
\begin{equation}\label{majks}
I(0)<   \mu^{**}
\end{equation}
is satisfied 
with 
\begin{equation}\label{muss}
 \mu^{**}=\mu^{**}(\alpha,M)=\frac{1}{\alpha} h_2(M/8\pi),
\end{equation}
where 
$h_2(s)=\pi.\frac{s(s-1)}{s+1}\left[ \log(2s/(s+1))\right]^2$, $s>1$,
then the solution $(n_t)$ blows up in finite time $T^*$ such that 
\begin{equation}\label{tks}
T^*\leq T^*_{ks}(M,I(0)):=\frac{I(0)}{M\left(\frac{M}{8\pi} -1\right)}.
\end{equation}
(See \cite[Theorem  2.9]{Ka-Su} for a generalization  of blow-up results to a larger class of kernels for  aggregation equations  including Bessel kernels).

Both bounds of blow-up criteria  \eqref{lamc} and \eqref{majks} on $I(0)$ are comparable when $M$
 is large enough up to a multiplicative constant.
But  note that $\mu^{**}(\alpha,M)<< \mu^*(\alpha,M)$ when $M$ is closed to $8\pi$.
The blow-up criteria $ \mu^*$  and $ \mu^{**}$ above are certainly not optimal when $M$ becomes large, as it   will   be shown below, 
in the sense that there exists a larger  $\mu$ with  $\mu>>\mu^*$ and $\mu>>\mu^{**}$  such that 
$T^*$ is finite. In fact, our main  result    will    increase the set  of initial data  $n_0$  for which the solutions of the (PKS) equation blow-up,
see Theorem \ref{bupp} below.

 \section{Main results}\label{mainresult}
\subsection{Aims of the paper}\label{aim}
In this paper, we  focus on the blow-up phenomenon for the Keller-Segel system on $\R^2$
with a  consumption term  $\alpha>0$ for the  supercritical mass case $M>8\pi$.
The aim of this study  is  to  show that we can enlarge the bounds of \eqref{ccond} on the second moment $I(0)$ of the initial data $n_0$
obtained in  \cite[Theorem 1.2]{CC1} and \cite[Theorem 2,p.356]{koz-sug} implying blow-up results.

More precisely, our aim is to find  largest moment bound  $\lambda^* (\alpha,M)$ such that $\mu^*(\alpha,M)$ and 
$ \mu^{**}(\alpha,M)$, given respectively by \eqref{lamc} and \eqref{muss}, 
satisfies
 $$
 \lim_{M\rightarrow +\infty} 
 \frac{\lambda^*(\alpha,M)}
 {\mu(\alpha,M)}=+\infty,
 $$
(also when ${M\rightarrow 8\pi^+}$) for each value
$\mu(\alpha,M)=\mu^*(\alpha,M)$ and $\mu(\alpha,M)= \mu^{**}(\alpha,M)$, and so that the condition
$I(0)<\lambda^* (\alpha,M)$ implies the blow-up.

 Moreover, we obtain  a precise  upper bound for the blow-up  time $T^*$ under this new condition 
 $I(0)<\lambda^* (\alpha,M)$. 
Our  constant $\lambda^* (\alpha,M)$    will   also take  into account of  the center of mass defined as follows,
\begin{equation}\label{bary}
 B_0=B_0(t)=\frac{1}{M} \int_{\Ri^2}   x.n_t(x)\, dx=\frac{1}{M}\int_{\Ri^2}   x. n_0(x)\, dx \in \R^2, \quad t\in (0,T^*).
\end{equation}
This is also the expected value  for  the density  of probability $n_0/M$ on $\R^2$. 
 The existence of the center of mass  is due to the assumption of  finiteness of the second moment 
 by  using  H\"{o}lder's inequality. 
The fact that the center $B_0(t)$ is constant with respect to the time $t$ is proved 
 by the usual method   of derivation of $B_0(t)$ and by using the (PKS) equation.

The center of mass $B_0$  arises in  the study  of  
 global dynamics of the  (PKS) equation on $\R^n$, $n\geq 3$ and $\alpha=0$, in  subcritical mass $M<8\pi$,
 see for instance the recent paper \cite[Theorem 1.4]{H-Y}. For the case $n=2$, the center of mass is not mentioned explicitly in Theorem 1.2 of \cite{H-Y}.
  Note  also that, in many papers, the center of mass  doesn't appear explicitly, or is fixed to be  $0\in \R^n$, e.g.  for radial solutions.
The center of mass has already  been  considered  in  the study of aggregation equations,
 see   for instance \cite{BCL}.
 
On the other hand,  the value of the mass $M$ of $n_0$ is  sometimes fixed to be $1$, i.e. the solution $(n_t)$ is a density of probability for each $t\geq 0$,
and, moreover the density is assume to be centered, i.e. $B_0=0$. This often simplifies the presentation of  theoretical computations.

It  is  certainly biologically relevant from the experimental point of view to take into account  
of the position  of the center of mass for the evolution of the cells,
but  also relevant  from a mathematical point of view, as we   will   see in Theorem \ref{bupp}. 
Of course,   it is natural  to impose the condition  $B_0=0$,  
that is  the center of mass  is at the origin of the plan,
in particular to simplify  numerical simulations issues, see \cite{tomth} for instance.  
\vskip0,3cm

Since the second moment and the center of mass of the initial data  $n_0$ arise in our discussion,
we are led to  consider  the variance $V_2(n_0)$ defined by
\begin{equation}\label{varzero}
V_2(n_0)=\frac{1}{M}\int_{\Ri^2} \vert x-B_0\vert^2 n_0(x)\, dx
\end{equation}
for the study of the (PKS) system.
Furthermore,  it seems  to be  more appropriate to use the variance
for the formulation of our results instead of the second moment $I(0)$  of $n_0$.
Indeed, there are   at least  two reasons for doing this.
The variance arises naturally in our method of  proof when using a convexity argument,
and also for the simplicity of the resulting statements.
The systematic  use of  variance $V_2(n_0)$ appears to be relatively  recent for the study of the (PKS) system,
to the authors know. Despite the fact that this notion, classic  in probability theory, underlies many works on the subject.
We   will  adopt this point of view in this paper.
\vskip0.3cm

Let us be more specific about the comparisons of our new results 
and previous results, expressed in terms of variance criteria.
For that purpose, the conditions of blow-up mentioned above 
and originally expressed with the second moment 
are now  reformulated in terms of variance, as follows. 
Let $M>8\pi$. We   will   assume below that $B_0=0$  to homogenize  the presentation of the results.

The blow-up criterion given by  \cite[Theorem 1.2]{CC1} takes  the following form,
\begin{equation}\label{gammacc}
V_2(n_0):=\frac{I(0)}{M}
\leq \gamma^*_{cc}(\alpha,M)
:=
\frac{\mu^*}{M}
=
\frac{(M-8\pi)^2}{4\alpha {C}^2 M^2},
\end{equation}
and  that of \cite[Th2.p.356]{koz-sug},
\begin{equation}\label{gammaks}
V_2(n_0)  
\leq \gamma^*_{ks}(\alpha,M)
:=
\frac{\mu^{**}}{M}
=\frac{1}{\alpha M} h_2(M/8\pi)
=
 \frac{1}{8\alpha}
 \frac{(M-8\pi)}{(M+8\pi)}
 \left[ \log\left (\frac{2M} {M+8\pi }\right)\right]^2.
\end{equation}

We can already  compare both  criteria above  with the  new  following blow-up  criterion 
\begin{equation} 
V_2(n_0) <  
\gamma^*(\alpha, M)=
\frac{1}{2\alpha}\left(
g_{1}^{-1}\right)^2
(8\pi /M),
\end{equation}
i.e. \eqref{eqgamma}, given  in Theorem \ref{bupp}.
Here $g_{1}^{-1}$ denotes the inverse function of $g_{1}$,  see  \eqref{galphafn} or  \eqref{defgalpha}.
\vskip0.3cm

Recall that the  blow-up criteria  \eqref{gammacc}-\eqref{gammaks}-\eqref{eqgamma} 
impose the smallness of the variance $V_2(n_0)$, i.e.  a large concentration of cells around the center of mass
at time $t=0$.  
Consequently, a blow-up in the evolution of the biological process is expected in finite time for supercritical mass $M>8\pi$ when
$\alpha >0$.
The larger  $\gamma^*$ is,   
the larger the set of initial data where a blow-up of the  corresponding solution of the (PKS) equation occurs.
 \vskip0.3cm

These criteria   can be compared in at least  two directions.

(i) When $M$ tends  to the critical mass $8\pi$ with  $M>8\pi$.

(ii) When $M$ tends  to $+\infty$.

Here are the details of  the comparison of the  ratios of  two  gammas chosen  from  the three values
$\gamma^*$, $\gamma_{cc}$ and $\gamma_{ks}$ defined above.

 Let $\alpha >0$ and $M>8\pi$ be fixed. 
\vskip0.3cm

(i) {\it Comparison when $M$  goes to $8\pi$}.

First, note that  both functions $\gamma^*_{cc}(\alpha,M)$, $\gamma^*_{ks}(\alpha,M)$   and 
$\gamma^*_{}(\alpha,M)$ tends to zero
when $M\rightarrow8\pi^+$ for fixed $\alpha>0$.
Moreover, we easily  have the following asymptotic for the ratio,
$$ 
\frac{\gamma_{cc}^*(\alpha,M)}{\gamma^*_{ks}(\alpha,M)}
\sim 
\frac{K}{M-8\pi}, \quad M\rightarrow 8\pi^+.
$$
Thus, we deduce that  
$$
\lim_{M\rightarrow 8\pi^+}  
\frac{\gamma_{cc}^*(\alpha,M)}{\gamma^*_{ks}(\alpha,M)}=+\infty.
$$
It  can  also be  shown that 
\begin{equation}\label{gamgamcc}
\lim_{M\rightarrow 8\pi^+}  
\frac{\gamma^*(\alpha,M)}{\gamma^*_{cc}(\alpha,M)}=+\infty.
\end{equation}
See Proposition \ref{progamcc}  in Section \ref{gestim} for a proof of this fact. (We also provide a more precise result, 
that is  a lower bound for the growth of this ratio as a function of $M$.)

Hence, this implies  that  ${\gamma^*(\alpha,M)}>> \gamma_{cc}^*(\alpha,M) >> \gamma_{ks}^*(\alpha,M)$,
as $M$ is  closed to $8\pi$. So, the blow-up criterion with  $\gamma_{cc}^*(\alpha,M)$ produces more blow-up situations than 
$\gamma_{ks}^*(\alpha,M)$, and  $\gamma_{}^*(\alpha,M)$ 
more blow-up situations than 
$\gamma_{cc}^*(\alpha,M)$,  when $M$ is close to $8\pi$.
\vskip0.3cm

(ii) {\it Comparison when $M$  goes to infinity}.

For large mass, i.e. as $M$ goes to infinity, we have the following asymptotics:
$$
{\gamma^*_{cc}(\alpha,M)}
\sim 
\frac{1}{4\alpha {\mathcal C}^2},
\qquad
{\gamma^*_{ks}(\alpha,M)}
\sim 
\frac{(\ln 2)^2}{8\alpha},
\qquad
{\gamma^*(\alpha,M)}
\sim 
\frac{(\ln M)^2}{2\alpha}.
$$
Note that the limit as $M\rightarrow +\infty$ of  $\gamma^*_{cc}(M)$ and $\gamma^*_{ks}(M)$ are both constants,  
unlike ${\gamma^*(\alpha,M)}$ which  is unbounded and which 
has logarithmic  squared  growth  as a function of $M$.
Thus,  the bound ${\gamma^*(\alpha,M)}$ differs not only quantitatively,
 but also qualitatively from  the bounds ${\gamma^*_{cc}(\alpha,M)}$ and 
${\gamma^*_{ks}(\alpha,M)}$  for large $M$.
For some  estimates of $g_1^{-1}(\rho)$ as $\rho$ goes to zero, corresponding to the study of the case where $M$ goes to infinity, see Section \ref{gestim}.
See also Corollary \ref{blowcreation} for an interesting consequence of this fact.
Hence, the  blow-up criterion $V_2(n_0) < \gamma^*_{}(\alpha,M)$
 leads to a larger set of initial data $n_0$ producing more  blow-up solutions of  the (PKS) equation than those obtained
 with ${\gamma^*_{cc}(\alpha,M)}$ or ${\gamma^*_{ks}(\alpha,M)}$.

 As  in  \cite{CC1} and \cite{koz-sug}, 
 we  provide an explicit  bound on the maximal existence time $T^*$ of
 solutions of (PKS) equations under the blow-up criterion $V_2(n_0) < \gamma^*_{}(\alpha,M)$. 
 The  value of $T^*$ is  bounded  by a function  $T^*_{\alpha}(n_0)$, 
 which decreases as the variance $V_2(n_0)$ decreases  as physically expected,
 see \eqref{Tboundzero}.

 \subsection{Main results on blow-up  time and  maximal existence time} \label{result}
We now  present the  main  result  of this paper, providing  a new criterion  
on the second moment (in fact, on the variance) of the  initial data $n_0$
for which the solutions  of the (PKS) equation blows up on  $\R^2$. 
This leads to more situations where the  blow-up occurs  than those obtained in the literature,
as far as  the authors knows.
Our approach takes into account the nature of the (PKS)  equation   in a much more  intrinsic way. 
The main argument takes advantage of the convexity of the function
  $\rho\rightarrow g_{\alpha}(\sqrt{\rho})$, where the function $g_{\alpha}$ is defined by  \eqref{defgalpha}  below. 
 \vskip0.3cm
 
Throughout  the paper, we assume that the mass $M$ satisfies the condition $M>8\pi$.
We    will  also    assume the local existence of solutions (in a certain sense) which allows us to study  the evolution of the 
the second moment of the solution $(n_t)$ for $0\leq t<T^*$,
assuming  the finiteness of the second moment of the  initial data $n_0$.

First, we recall that  the  function $g_{\alpha}$ that plays a fundamental  role in our discussion,
\begin{equation}\label{defgalpha}
g_{\alpha}(r)=\int_0^{+\infty}e^{-\frac{\alpha r^2}{4s}}  e^{-s}\, ds, \quad r\geq 0, \quad \alpha \geq 0.
\end{equation}
The importance of this function is given by  the relationship \eqref{linkbesselgalpha} between  
the gradient of the Bessel kernel $\nabla B_{\alpha}$ and  $g_{\alpha}$.
It is clear that  the function $\rho\rightarrow g_{\alpha}(\sqrt{\rho})$ is   convex, which is  crucial 
for  the main result, i.e.  Theorem \ref{bupp} of this paper.
Note also that the  function $g_{\alpha}$ is strictly decreasing  and its range is 
$(0,1]$ for all $\alpha>0 $. 
Hence, the inverse  function $g_{\alpha}^{-1}:[0,1]\rightarrow [0,+\infty)$ is well-defined, strictly decreasing
and $\lim_{u\rightarrow 0} g_{\alpha}^{-1}(u)=+\infty$.
We will   sometimes consider the case $\alpha=0$, and the corresponding function
$g_{0}\equiv 1$ on $[0,+\infty)$ as a limiting case.
\vskip0.3cm

The main result of this paper is the next  one.

\begin{theo}\label{bupp}
Let $\alpha>0$.
Assume that $(n_t)_{0<t<T^*}$ is a solution of the (PKS) equation  \eqref{ks} with  non-negative initial  data
$n_0\in L^1(\R^2)$ of mass $M=\int _{\Ri^2} n_0(x)\, dx>8\pi$.
Assume that  the second moment $I(0)$ of $n_0$, i.e.  \eqref{secondmoment}, is finite, and 
that the variance $V_2(n_0)$  \eqref{varzero} of $n_0$ satisfies the next  inequality,
\begin{equation}\label{eqgamma}
V_2(n_0)
<
\gamma^*(\alpha, M):=
\frac{1}{2\alpha}\left(
g_{1}^{-1}\right)^2
(8\pi /M),
\end{equation}
where $g_{1}^{-1}$ is the inverse function of $g_{\alpha}$ for  $\alpha=1$ (see \eqref{defgalpha}).
\begin{enumerate}
\item
Then the solution $(n_t)$ of the (PKS) equation blows up in finite time, i.e. $T^*<+\infty$,
where $T^*$ is the maximal  existence time  of the solution. 
Moreover,  we  have the following explicit bound  on $T^*$,
\begin{equation}\label{Tboundzero}
T^*\leq T^*_{\alpha} (n_0):=
2\pi
\int_0^{V_2(n_0)} \frac{ds}
{
M g_{1}\left(\sqrt{2\alpha s}\right)
-8\pi
}
<+\infty.
\end{equation}
\item
 In particular,  the maximal existence time $T^*$  satisfies
\begin{equation}\label{Tbound}
T^*\leq 
T^*_{\alpha} (n_0)
\leq 
\frac{2\pi V_2(n_0)}
{
M g_{1}\left(\sqrt{2\alpha  V_2(n_0)}\right)
-8\pi
}
.
\end{equation}

\item
Moreover, we  have the following bound on the time  evolution of the variance $V(t)$ of the solution $(n_t)$,
\begin{equation}\label{vtheta}
V(t):=
\frac{1}{M}\int_{\Ri^2} \vert x-B_0\vert^2 n_t(x)\, dx
\leq
\Theta^{-1}(t), \quad 0<t<T^*,
\end{equation}
where $\Theta^{-1}$ is the inverse function of
$$
\Theta(x)
=
2\pi
\int_x^{V_2(n_0)} \frac{ds}
{
M g_{1}\left(\sqrt{2\alpha s}\right)
-8\pi
}
,\quad  0\leq x \leq V_2(n_0),
$$
and the derivative of the variance  is controlled as follows,
\begin{equation}\label{vprime}
V^{\prime}(t)
\leq
f(\Theta^{-1}(t)), \quad 0<t<T^*,
\end{equation}
where 
$$
f(\lambda)
=
4-
\frac{M}{2\pi}g_1\left(\sqrt{2{\alpha}\lambda}\right),
 \; \lambda>0.
$$
\item
In particular,  we have  for all $ 0\leq t<T^*$,
\begin{equation}\label{weakboundks}
V (t)\leq   V(0)+tf(V(0))
\quad {\rm and} \quad 
V^{\prime}(t)\leq f(\, V(0)+tf(V(0)) \,),
\end{equation}
where $V(0):=V_2(n_0)$ and $f(V(0))<0$.
\end{enumerate}
\end{theo}

Before proving  Theorem \ref{bupp} in Section \ref{proofth}, we make several remarks and comments.

\begin{rem}
Note that $V^{\prime}(t)\leq 4$ for all $t\in (0, T^*)$.
\vskip0.3cm
\end{rem}

\begin{rem}
The blow-up condition \eqref{eqgamma} on the variance $V_2(n_0)$ is equivalent to 
the following but less elegant  condition on the second moment,
\begin{equation}\label{izerocond}
I(0)=
\int_{\Ri^2}  \vert x\vert^2\,n_0(x)\, dx
\leq
\frac{M}{2\alpha}\left(
g_{1}^{-1}\right)^2
(8\pi /M)+ M\vert B_0\vert^2.
\end{equation}
\vskip0.3cm
\end{rem}

\begin{rem}
Instead of studying  the variation of the second moment $I(t)$ given by \eqref{ioft}, 
we prefer to study 
the time evolution    of the variance  given by
$$
(0,T^*) \ni t \mapsto V(t)=\int_{\Ri^2}\vert x- B_0\vert^2 n_t(x)\,\frac{dx}{M}
=\frac{I(t)}{M}-\vert B_0\vert^2.
$$
So, we have 
\begin{equation}\label{varmom}
V_2(n_0)=\frac{1}{M} \int_{\Ri^2} \vert x-B_0\vert^2 n_0(x)\, dx
=\frac{I(0)}{M}-  \vert B_0\vert^2,
\end{equation}
where $B_0=\int_{\Ri^2}  x\,n_t(x)\, dx/M$ is the center of mass  of $n_0$.
\vskip0.3cm
In the  framework of weak solutions for instance, it is proved in   \cite[Cor.2.2]{BDP},
 under finite second moment  and finite entropy assumptions, and when  $\alpha=0$
and $M>8\pi$, that  the maximal existence time $T ^*$ of the solutions
 is bounded as follows,
\begin{equation}\label{bupbdp}
 T^*\leq 
 T^*_m:=\frac{2\pi  I(0) }{M(M-8\pi)} < + \infty,
\end{equation}
It  can  easily be shown that  \eqref{bupbdp} also holds true with $ \frac{I(0)}{M}$ 
replaced in the expression of  $T^*_m$  by the a priori smaller quantity $V(0)=V_2(n_0)$,  
i.e.
\begin{equation}\label{bupbdp2}
 T^*\leq T^*_v:= \frac{2\pi V(0)}{(M-8\pi)}.
\end{equation}
So, in this article  and for  the general situation, i.e.  $\alpha\geq 0$,  it seems 
more appropriate to  also include  the center of mass $B_0$ in the discussion.
Indeed, if $B_0\neq 0\in \R^2$ then we have  $V_2(n_0):=V(0)<\frac{I(0)}{M}$, 
and consequently  $T^*_v <T^*_m$, which improves 
the control of $ T^*$ in  \eqref{bupbdp}. 
Note that the difference  $T^*_v - T^*_m$  of time bounds  of $T^*$  
can be simply quantified as follows,
 $$ 
 T^*_v-T^*_m= :\frac{- 2\pi \vert B_0\vert^2}{(M-8\pi)}.
 $$
The bound   $T^*_v$  above on the maximal existence time $T^*$ 
seems to be  rarely mentioned  to the best  of our knowledge (even for the case $\alpha= 0$).
 \vskip0.3cm
\end{rem}

\begin{rem}
 All the bounds  on the maximal existence time $T^*$ mentioned in that paper, i.e 
$T^*_{cc}$, $T^{*}_{ks}$ and $T^*_{\alpha}$
 given respectively by \eqref{tcc},\eqref{tks} 
and  \eqref{Tboundzero},  lead  formally  at the  limiting  case  $\alpha\rightarrow  0^+$,
to the well-known bound on   $T^*$ for  $\alpha=0$, i.e. \eqref{bupbdp} above 
when the center of mass  $B_0$ satisfies  $B_0=0\in \R^2$.
This shows  the asymptotic  sharpness of all these estimates when $\alpha$ is closed to $0$ with respect to the usual bound  \eqref{bupbdp}.
Note  that the upper  bound assumptions   \eqref{gammacc}-\eqref{gammaks}-\eqref{eqgamma} on $I(0)$ then 
vanish when  formally $\alpha=0$.
 \vskip0.3cm
\end{rem}

\begin{rem}
For the convenience of the reader, we recall the arguments of the proof  of  statements \eqref{bupbdp} and \eqref{bupbdp2} when $\alpha=0$.
We  begin with the following well-known formula,
 $$
 I^{\prime}(t)=\frac{d}{dt}\left( \int_{\Ri^2} \vert x\vert^2 n_t(x)\, dx\right)
 =4M\left(1-\frac{M}{8\pi}\right),
 $$
which holds  for  any  solution $(n_t)$  of  the (PKS) equation with  a finite second moment assumption on $n_0$.
 This equality can also  be written in terms of the  variance  as follows,
 $$
 V^{\prime}(t)=\frac{1}{M}  I^{\prime}(t)
 =4\left(1-\frac{M}{8\pi}\right),
 $$
where $V(t)=\frac{1}{M}\int_{\Ri^2} \vert x-B_0\vert^2 n_t(x)\, dx$.
After integration with respect to the variable $t$, 
we get 
 $$
 V(t)=V(0)+4\left(1-\frac{M}{8\pi}\right)t, \quad 0<t<T^*.
 $$
Since $V(t)$  is  non-negative for  all $0<t<T^*$,  and $M>8\pi$,
this  implies that $T^*$ is finite and  \eqref{bupbdp2} holds true.
(See \cite[p.2808]{FJ} where this result can be found).
As a consequence, we obtain \eqref{bupbdp},  since
$V(0)=\frac{I(0)}{M}-  \vert B_0\vert^2\leq \frac{I(0)}{M}$.
\end{rem}
\vskip0.3cm

\begin{rem}
Note that  the variance $V(0):=V_2(n_0)$ can be taken  as  close to zero  as we want by  an appropriate choice of $n_0$.
But this cannot be the case  for the second moment $I(0)$ if $B_0\neq 0$.
Indeed,  we  have  
\begin{equation}\label{lower}
I(0)\geq M\vert B_0\vert^2,
\end{equation}
since $0\leq V(0)=\frac{I(0)}{M}-  \vert B_0\vert^2$. 
In particular, the explicit bound  $T^*_{\alpha}(n_0)$ on  $T^*$, given by \eqref{Tboundzero}, says that if the variance is small,
 i.e. the cells are   distributed near  the center of mass, then the blow-up time $T^*$ is small.  
Of course, this phenomenon  is biologically expected in the case of supercritical mass $M>8\pi$, 
and it can be  quantified by \eqref{Tboundzero}.
Note that this fact cannot be  directly observed from  inequality \eqref{bupbdp} due to the lower bound \eqref{lower} on  $I(0)$ described above.
\end{rem}

\begin{rem}
For $\alpha>0$ and for the blow-up issue,  note that the variance $V(0)$, or the moment $I(0)$) condition \eqref{eqgamma},
 is now directly  related to the (PKS)  equation 
\eqref{ks} through the Bessel kernel $B_{\alpha}$ via  $g_{\alpha}$. 
For the explicit relationship between $I^{\prime}(t)$ and $g_{\alpha}$, see  \eqref{igalpha}.
\end{rem}

\begin{rem}
Theorem  \ref{bupp} improves substantially the results of \cite[Theorem 1.2]{CC1},  and \cite[Theorem 2,p.356]{koz-sug}
when   $M$ is large, and  also when $M$ is closed to $8\pi$, providing more situations of blow-up.
See Section \ref{aim} above for a precise discussion.
\vskip0.3cm
\end{rem}

\begin{rem}
For all given $\alpha>0$ and all $M>8\pi$,  there exists initial conditions $n_0$ satisfying \eqref{eqgamma}, for instance 
$n_0=c\chi_{B_{\varepsilon}}$ for ${\varepsilon}$ small enough  with $c>0$ large enough,
 such that the mass $\int n_0\,dx= M>8\pi$. Here $\chi_{B_{\rho}}$ denotes the characteristic function 
of the ball $B_{\rho}$ of  radius ${\rho}$ and center  (of mass) $B_0=0\in \R^2$.
\vskip0.3cm
\end{rem}

\begin{rem}
Let  $n_0$ be  any nonnegative initial data  fixed with finite second moment and mass  $M>8\pi$.
Then, for any $\alpha$ such that
\begin{equation}\label{interv}
\alpha \in \left(0, 
\frac{1}{2V_2(n_0)}\left(
g_{1}^{-1}\right)^2
(8\pi /M)\right),
\end{equation}
the condition  \eqref{eqgamma} is fulfilled. Hence,  the  corresponding solution  of the (PKS)   equation  blows up in finite time.
(Recall that the function $g_{1}$  is defined by formula  \eqref{defgalpha}.)
\vskip0.3cm
It is proved in \cite[Theorem 2.2]{BCKZ} that for any nonnegative  initial data $n_0$ (with $M>8\pi$, or not),
 there exists a global-in-time  (mild) solution of the $(PKS)$ equation
 if $\alpha\geq \alpha(n_0)$ with $\alpha(n_0)$ large enough.
No finite second moment condition has to be imposed on $n_0$.
Of course, there  is no contradiction with our  blow-up condition   \eqref{eqgamma}.
In the case of finite  second moment for $n_0$, the condition  \eqref{interv} 
provides an  explicit interval for  $\alpha$ on   which the  blow-up phenomenon occurs. 
Hence, with the additional   condition of finite second moment on the initial data, we deduce  in the situation described  above  that  
$$
 \alpha(n_0)\geq 
\frac{1}{2V_2(n_0)}\left(
g_{1}^{-1}\right)^2
(8\pi /M).
$$
\vskip0.3cm
\end{rem}

\begin{rem}
In \cite[Theorem 2.2 (ii)]{BCKZ}, {local} criteria are given on nonnegative initial data $n_0$ (without a priori   finite second moment assumption on $n_0$)
leading to a blow-up of the corresponding solution of the (PKS) equation on $\R^2$. 
The authors use a different approach with respect to 
{global} criteria  usually considered, for instance in   \cite{CC1}, \cite{koz-sug}, and also  in the present paper.
\vskip0.3cm
\end{rem}

\begin{rem}
The quantity $V(0)=V_2(n_0)$ can  be interpreted as  the variance of a 2-dimensional random variable 
$X_0$ with density  law  $n_0/M$ with respect to Lebesgue measure. 
We then  obtain a blow-up when  we impose a control  from above on the  variance of $X_0$,
that is to say  a condition on  the spreading  (not too large)  of the cells at the initial time $t=0$.
The smaller $V_0$,  the more  the cells   are concentrated around the center of mass  
$B_0 $ at the beginning of the evolution process.
For a probabilistic interpretation of the (PKS)  system, see for instance  \cite{FTom} and references therein. 
See also \cite{FT} for the subcritical case.
\vskip0.2cm
\end{rem}

\begin{rem}
The condition of blow-up \eqref{eqgamma}  could be sharp,
but we have no evidence of this fact.
\end{rem}

{\bf{Proof of Theorem \ref{bupp}}\label{proofth}}
In this proof, we  assume that $M>8\pi$ and  $\alpha >0$.
\vskip0.2cm
(1)  
{\em Step 1}.
Let  $I(t)$ be  the second moment of the solution $n_t$ to the (PKS) equation   defined by
$$
I(t)=\int_{\Ri^2} \vert x\vert^2 n_t(x)\, dx, \; 0<t<T^*.
$$
From \eqref{ks} and \eqref{linkbesselgalpha},  we  successively deduce that  the derivative of  $I$ satisfies
the next equation
$$
\frac{d}{dt} I(t)=I^{\prime}(t)=4M+2 \int_{\Ri^2} n_t(x)  x.(\nabla B_{\alpha}\ast  n_t)(x)
\,dx.
$$
Calculating  the gradient of $B_{\alpha}$ defined in \eqref{besselkern}, this leads us to the next expression
$$
I^{\prime}(t)=4M-\frac{1}{\pi}\int_{\Ri^2}\int_{\Ri^2} 
 n_t(x) 
\frac{x.(x-y)}{\vert x-y\vert^2}
g_{\alpha}(\vert x-y\vert) n_t(y)\,dxdy.
$$
By symmetry in $(x,y)$ in the double integral, we get
\begin{equation}\label{igalpha}
I^{\prime}(t)=4M-\frac{1}{2\pi}\int  \int g_{\alpha}(\vert x-y\vert) n_t(x)n_t(y)\,dxdy.
\end{equation}
 (Our notation is slightly different from  \cite{CC1} where our $ g_{\alpha}(\vert x-y\vert)$ 
corresponds to $ g_{\alpha}(x-y)$  of  \cite[Section 7]{CC1}).
 \vskip0.3cm

 {\em Step 2}.
 For any $t>0$ , we define the measure $d\nu_t$ by
$$
d\nu_t(x,y)=n_t(x)n_t(y)\,dxdy/M^2.
$$
This measure is a probability measure on the product space $\R^2\times \R^2$.
For any $\rho\geq 0$, we set 
$$
\Psi_{\alpha}(\rho)=
-\int_0^{+\infty} e^{-s} e^{-\frac{\alpha \rho}{4s}}\,ds
=
-g_{\alpha}(\sqrt{\rho}).
$$

Thus, we can rewrite the derivative $I^{\prime}(t)$ as follows
$$
I^{\prime}(t)=4M+\frac{M^2}{2\pi} 
\int\int
\Psi_{\alpha}(\vert x-y\vert^2)d\nu_t(x,y).
$$
Since $d\nu_t$ is a probability measure and $\Psi_{\alpha}$ is a concave (increasing) function 
on $(0,+\infty)$,  we can apply Jensen's inequality  which asserts  that
$$
 \int\int 
\Psi_{\alpha}(\vert x-y\vert^2)\,d\nu_t(x,y)
\leq
\Psi_{\alpha}\left(  \int \int \vert x-y\vert^2\,d\nu_t(x,y)\right).
$$
So, we deduce that
$$
I^{\prime}(t)
\leq 
4M+\frac{M^2}{2\pi} 
\Psi_{\alpha}\left(  \int \int \vert x-y\vert^2\,d\nu_t(x,y)\right).
$$
Next, we use the following formula
$$
\int  \int \vert x-y\vert^2\,d\nu_t(x,y)= \frac{2}{M}(I(t)-M\vert B_0\vert^2),
$$
where  $B_0=\int x\,n_0(x)\, dx/M$.
(In particular, this yields a lower bound on $I(t)$, i.e.
$I(t)\geq M\vert B_0\vert^2$.
Thus $B_0$ is finite since $M\vert B_0\vert^2\leq I(0)<\infty$.)

Then the differential inequality just above can be read  as follows
\begin{equation}\label{diffit1}
I^{\prime}(t)\leq
4M+\frac{M^2}{2\pi}
\Psi_{\alpha}\left( \frac{2}{M}(I(t)-M\vert B_0\vert^2)\right)
\end{equation}
for all $0<t< T^*$.
\vskip0.3cm

{\em Step 3}.
Let  $V(t)$  be the variance at time $t\in (0,T)^*$  of $(n_t)$ defined by
\begin{equation}\label{vartimet}
V(t)=\frac{1}{M}\int_{\Ri^2} \vert x-B_0\vert^2 n_t(x)\, dx.
\end{equation}
For all  $0\leq t<T^*$, we have 
$$
V(t)= \frac{1}{M}(I(t)-M\vert B_0\vert^2).
$$
With our notation above, we can write  $V_2(n_0)=V(0)$ (see \eqref{varzero}).
Recall that we also have  for the center of mass $B(t)=B(0)$ for all  $0\leq t<T^*$.
In terms of $V(t)$, the differential equation \eqref{diffit1} can be written as
\begin{equation}\label{diffit2}
V^{\prime}(t)\leq
4+\frac{M}{2\pi}
\Psi_{\alpha}\left( 2V(t)\right),
\end{equation}
for all $0<t<  T^*$.
For all  $\lambda>0$, we define  $f$  by the following expression 
$$
f(\lambda)=
4+\frac{M}{2\pi}
\Psi_{\alpha}\left(2 \lambda \right)
=
4-
\frac{M}{2\pi}g_{\alpha}\left(\sqrt{2\lambda}\right),
$$
or more explicitly,  
\begin{equation}\label{fp}
f(\lambda)=
 4-
\frac{M}{2\pi} \left(
\int_0^{+\infty}
e^{-s} e^{-\frac{\alpha \lambda }{2s}}\,ds
\right).
\end{equation}
This function $f$ depends only on $M$ and $\alpha$.
The differential inequality just above can be written in a final form as,
 \begin{equation}\label{vprime2}
V^{\prime}(t) 
\leq
f(V(t)), \quad 0<t<T^*,
\end{equation}
where $V(t)$  is defined by  \eqref{vartimet}.
\vskip0.3cm

{\em Step 4}.
From \eqref{vprime2}, we are now in position to apply Proposition \ref{diffineq} 
and conclude the proof of the theorem.
The function  $f:(0,+\infty)\rightarrow \R$ defined above  is  a strictly increasing continuous (concave)  function with
$f(0^+)<0<f(+\infty)=4\leq +\infty$. Note that the condition $f(0^+)<0$ is precisely the assumption  $M>8\pi$ (for any $\alpha >0$).

Let $\lambda^*>0$ be the unique zero of $f$, i.e. $f(\lambda^*)=0$.
Then the  value of $\lambda^*$ is given explicitly  by
\begin{equation}\label{lambdastar}
\gamma^*(\alpha, M):=\lambda^*= \frac{1}{2\alpha}\left(g_{1}^{-1}\right)^2(8\pi /M).
\end{equation}
Indeed, here we have used the fact that the function  $r\mapsto g_{\alpha}(r)$  from $[0,+\infty)$ to $(0,1]$ is strictly decreasing,
hence invertible. It also  satisfies the relation
$g_{\alpha}(r)=g_{1}(\sqrt{\alpha} r)$ for all $r\geq 0$ and all $\alpha>0$.
In particular, we have for all $y\in (0,1]$,
$$
g_{\alpha}^{-1}(y)=
\frac{1}{\sqrt{\alpha}}g_{1}^{-1}(y).
$$
The differential inequality \eqref{ineqdiff} of  Proposition \ref{diffineq}  
is satisfied by  $V(t)$, $0<t<T^*$ with $f$ given above. 
So, let $\Theta$ be defined by
$$
\Theta(x)=\int_x^{V(0)} \frac{ds}{-f(s)}
=
2\pi
\int_x^{V(0)} \frac{ds}
{
M g_{\alpha}\left(\sqrt{2s}\right)
-8\pi
}
, \quad  0\leq x\leq V(0).
$$
By \eqref{tcgeneral}, the  maximal existence time $T^*$  of the solution to the (PKS) equation  is bounded as follows,
 \begin{equation}\label{}
T^*\leq T^*_c:=\Theta(0)
=2\pi
\int_0^{V_2(n_0)} \frac{ds}
{
M g_{1}\left(\sqrt{2\alpha s}\right)
-8\pi
}
<+\infty.
\end{equation}
Thus Statement  (1)  and, in particular,  inequality \eqref{Tboundzero} of Theorem \ref{bupp} are proved.
\vskip0.3cm

By inequality \eqref{boundTV(0)} of Proposition \ref{diffineq} and $V(0)=V_2(n_0) $, we also obtain
$$
T^*\leq T_c^{**}
\leq 
\frac{V(0)}{-f(V(0))}
=\frac{2\pi V(0)}
{Mg_{\alpha}(\sqrt{2V(0)})-8\pi
}
$$
$$
=
2\pi V_2(n_0)
\left[
M
\int_0^{+\infty}
\exp\left(
-\frac{\alpha V_2(n_0)}
{2s }
\right)
e^{-s} 
\,ds 
-8\pi
\right]^{-1}.
$$
This  proves \eqref{Tbound} and conclude the proof of Statement (2) of Theorem \ref{bupp}.
\vskip0.3cm

(3-4)
Inequalities  \eqref{vtheta}-\eqref{vprime}-\eqref{weakboundks}  follow respectively  from \eqref{thetaboundgeneral}
and  \eqref{weakbound} of Proposition \ref{diffineq}.
The proof of Theorem \ref{bupp} is now completed.
\hfill $\square$
\vskip0.3cm

\subsection{Blow-up  time bounds and bounds on solutions }\label{diff-argument}
To prove  blow-up for a PDE's system, one  often uses blow-up  time   bounds obtained from the study 
of a one-parameter differential inequality  associated with  the system.
In most articles, the arguments are given during the proof. They are more or less specific to the differential inequality
under study.  

In this section, we  state a general  sharp result of blow-up  time   bounds for a  certain type of one-parameter differential inequalities.
The statement is formulated independently of the PDE's systems  from which they may arise.
These estimates are of independent interest and can be used for further references.
Along with obtaining  precise   blow-up  time   bounds,
 we also  provide  upper bounds on the function (and its first derivative) satisfying the   differential inequality.
Our result  below is therefore stated in an independent proposition,  to which we   will    refer in this article.

\begin{pro}\label{diffineq}
Let $V:[0,T^*)\rightarrow[0,+\infty)$ be a non-negative continuous function which is a  $C^1$-function  on $(0,T^*)$, 
with (maximal) domain $[0,T^*)$ where $T^*\leq +\infty$.
Let  $f:(0,+\infty)\rightarrow \R$ be a strictly increasing continuous function such that
$$
-\infty\leq f(0^+)<0<f(+\infty)\leq +\infty.
$$
Let  $\lambda^*>0$ be the unique zero of $f$, i.e. $f(\lambda^*)=0$.
Assume that the following inequalities are satisfied, for all $0<t<T^*$,

\begin{equation}\label{ineqdiff}
V^{\, \prime}(t)\leq f(V(t)) \quad  {\rm and} \quad  V(0)<\lambda^*.
\end{equation}
Let $\Theta$ be defined by
$$
\Theta(x)=\int_x^{V(0)} \frac{ds}{-f(s)}, \quad  0\leq x\leq V(0).
$$

\begin{enumerate}
\item
 Then the maximal  existence time $T^*$ is finite  and satisfies 
 \begin{equation}\label{tcgeneral}
T^*\leq T^*_c:=\Theta(0)
=
\int_{0}^{V(0)} \frac{ds}{-f(s)}
<+\infty.
\end{equation}
Moreover, we have 
\begin{equation}\label{thetaboundgeneral}
V (t)\leq \Theta^{-1}(t), \quad V^{\prime}(t)\leq f( \Theta^{-1}(t)),
\end{equation}
for all $0<t<T^*$. Here $\Theta^{-1}$ denotes the inverse function of $\Theta$.
\item
In particular, the maximal  existence time $T^*$  also satisfies
\begin{equation}\label{boundTV(0)}
T^*\leq T_c^{**}:=\frac{V(0)}{-f(V(0))},
\end{equation}
and
\begin{equation}\label{weakbound}
V (t)\leq   V(0)+tf(V(0))
, \quad V^{\prime}(t)\leq f( V(0)+tf(V(0))), 
\end{equation}
for all $ 0\leq t<T^*$.
\end{enumerate}
\end{pro} 

For the sake of completeness, we give  a  detailed proof of this proposition. 
 In particular, we avoid some of the usual arbitrary and unnecessary choices and/or approximations
  in the proof in order to obtain optimal estimates. 
 As a result, we improve on the usual estimates of such bounds for applications.
\vskip0.3cm

{\bf Proof}. Assume that $V(0)<\lambda^*$. 
\vskip0.3cm

(1) Our first  goal is to prove  that   we have $f(V(t))<0$ for all $ 0<t<T^*$ whenever $V(0)<\lambda^*$.
\vskip0.3cm
{\bf (a)}
{\bf (i)}
Let consider the set defined by
$$
{\mathcal F}=\{ t^{\prime} \in (0,T^*]: V(t)\leq \lambda^*, 0\leq t<t^{\prime}   \}.
$$
Let $T_0:=\sup{\mathcal F}$.
We   will    prove below that $T^*=T_0$ and $T_0\in {\mathcal F}$.
Assume   for a while that this result holds true, i.e. 
$$
 V(t)\leq \lambda^*, \quad  0\leq t<T^*.
 $$
Our assertion  $f(V(t))<0$ for all $ 0<t<T^*$,    will    be deduced as follows.
From   the first  assumption of  \eqref{ineqdiff}, we obtain
$$
V^{\prime}(t)\leq f(V(t))\leq f(\lambda^*)=0, \quad 0\leq t<T^*,
$$ 
because $f$ is non decreasing, $ V(t)\leq \lambda^*, 0\leq t<T^*$, and $f( \lambda^*)=0$. So, the function $V$ is non increasing on $[0,T^*)$.   
In particular,  we deduce that  $V(t)\leq V(0)<\lambda^*$ for all  $0\leq t<T^*$.  
Since  $f$ is strictly increasing,  we obtain 
\begin{equation}\label{ineqnegativ}
 f(V(t))\leq f(V(0))
< 0=f(\lambda^*),
\end{equation}
 for all  $0\leq t<T^*$, which is the expected assertion for (1) above.
\vskip0.3cm

Next, we prove the assertions  concerning ${\mathcal F}$ defined  above.
\vskip0.3cm

 {\bf (ii)} 
 We first show that the set ${\mathcal F}$ is not empty.
 Indeed, under the assumption $V(0)<\lambda^*$, and  by continuity of the function $V$ at $0^+$,
  there exists $t_1\in (0,T^*]$ such that 
$V(t)\leq \lambda^*$ for all $0\leq t<t_1$. 
Thus $t_1\in {\mathcal F}$ and $ {\mathcal F}$ is a non empty set. 
\vskip0.3cm

Now, consider  $T_0=\sup\{t', t'\in {\mathcal F}\}$, which  exists  in $]-\infty,+\infty]$
since ${\mathcal F}$ is not empty.
More precisely,  $T_0$  belongs to  the set $[t_1,T^*] \subset (0,T^*]$  since ${\mathcal F}\subset  (0,T^*]$ and $t_1\in {\mathcal F}$.
Next, we prove that  $T_0$ finally belongs to ${\mathcal F}$.
In other words, we have  to show that  $V(t)\leq \lambda^*$ for all $0\leq t<T_0$.

Let $(s_n)_n$ be an increasing  sequence of ${\mathcal F}$ converging to $T_0$, and let $0\leq t<T_0$. 
By convergence of the   sequence $(s_n)_n$  to $T_0$, there exists $s_{n_0}$ 
 such that $0\leq t<s_{n_0}$.
Because $s_{n_0}\in {\mathcal F}$ and $t\in [0,s_n)$, we  deduce that
$V(t) \leq \lambda^*$, and this proves that  $T_0\in {\mathcal F}$. 
\vskip0.3cm
  
{\bf (iii)}
 Next, we show that $T_0=T^*$.
On the contrary, assume that $T_0<T^*$.
In particular, the value  $V(T_0)$ is well defined.
By  assumption \eqref{ineqdiff} and the fact that $T_0\in {\mathcal F}$,  this implies that
$$
V^{\prime}(t)\leq f(V(t))\leq f(\lambda^*)=0, \quad 0\leq  t<T_0,
$$
since $f$ is non decreasing.
Hence,  the function $V$ is non increasing on $[0,T_0)$. 
We deduce  successively that  $V(t)\leq V(0)$ for all $0 \leq t< T_0<T^*$, and by continuity of $V$ at $T_0<T^*$, we deduce that  
$$
V(T_0)=\lim_{t\rightarrow T_0^-} V(t)
\leq V(0)<\lambda^*.
$$ 
Thus,  we have  $V(T_0) <\lambda^*$. By continuity of $V$ at $T_0<T^*$, there exists $T_1$ such that $0<T_0<T_1\leq T^*$,
and $V(t)\leq \lambda^*$ for all $T_0\leq t<T_1$. Hence, the inequality  $V(t)\leq \lambda^*$  
also  holds for $t\in [0,T_1)$, i.e.   $T_1\in {\mathcal F}$ with $T_0<T_1$.
This contradicts the maximality of $T_0$ in  ${\mathcal F}$. 
We finally conclude that $T_0:=\sup {\mathcal F}= T^* \in {\mathcal F}$. 
\vskip0.5cm

{\bf (b)} We now prove the second part of (1) of Proposition \ref{diffineq}.
Recall that we  assume that  $V(0)<\lambda^*$. We have  just proved in (a) that, for all $s\in [0,T^*)$,
$$
V^{\,\prime}(s)\leq f(V(s)) <0.
$$
This  implies that
$$
 \frac{\vert V^{\prime}(s)\vert }{-f(V(s))}
 =
  \frac{V^{\,\prime}(s)}{f(V(s))}\geq 1.
 $$
By integration over $s$, we obtain  the next inequality for any  $t\in [0,T^*)$,
\begin{equation}
\int_0^t \frac{ \vert V^{\prime}(s)\vert }{-f(V(s))}\, ds \geq t.
\end{equation}
The function $V$  is strictly decreasing since $ V^{\prime}<0$ on  $[0,T^*)$.
Hence, the function $V$ is a bijection from $[0,t]$ to $[V(t), V(0)]$ for any  $t\in [0,T^*)$.
By the change of variables $y=V(s)$, and  from the definition of 
$\Theta$ given  in Proposition \ref{diffineq},  we get 
the next formula 
 \begin{equation}\label{ineqtheta}
\Theta(V(t)):=
\int_{V(t)}^{V(0)} \frac{1}{-f(y)}\, dy
=
\int_0^t \frac{ \vert V^{\prime}(s)\vert }{-f(V(s))}\, ds
 \geq t.
\end{equation}
Now, we are in position to conclude  part (1) of  the proposition.

(i) Since $V\geq 0$ and $-1/f(y)>0$ for $y\in [0,V(0)]$, we have
$$
t\leq \Theta(V(t))\leq \int_{0}^{V(0)} \frac{-1}{f(y)}\, dy=\Theta(0),
$$
for all $0<t<T^*$. Hence, by letting $t\rightarrow T^*$, we obtain  \eqref{tcgeneral}, i.e $T^*\leq \Theta(0)$. 
The value of the integral $\Theta(0)$ is finite because the function $y\mapsto -f(y)$ is continuous and strictly  positive on the compact set 
$[0,V(0)]$.
\vskip0.3cm

(ii)
Clearly, the function  $\Theta:[0, V(0)]\longrightarrow [0,\Theta(0)]$
is given by
$$
\Theta(x)= \int_{x}^{V(0)} \frac{-1}{f(y)}\, dy, \; 0\leq x\leq V(0),
$$
 is a strictly decreasing  bijection
due to the fact that $\Theta^{\prime}=1/f<0$ on $(0, V(0))$.
Hence, $\Theta$ is invertible with inverse $\Theta^{-1}$.
From inequality  \eqref{ineqtheta}, we first get that
$$
V (t)
\leq \Theta^{-1}(t),\quad 0<t<T^*.
$$
By applying  the assumption \eqref{ineqdiff}, we also deduce that
$$
 V^{\prime}(t)\leq f( \Theta^{-1}(t)), \quad 0<t<T^*,
 $$
 since $f$ is a non decreasing function.
So, we have proved   inequalities \eqref{tcgeneral} and \eqref{thetaboundgeneral}.
This concludes part (1) of  the proposition.
\vskip0.3cm

(2) This part is a direct and an  easy  consequence of part (1).
Indeed, from \eqref{tcgeneral}   
and the fact that $-1/f$ is non decreasing, we  get
$$
T^*\leq  T^*_c:=\Theta(0)= \int_{0}^{V(0)} \frac{1}{-f(y)}\, dy
 <
 T_c^{**}:=\frac{V(0)}{- f(V(0))}.
 $$
This gives  the bound \eqref{boundTV(0)}
on the maximal existence time  $T^*$  (and again the  finiteness of $T^*$).
On the other hand,  from \eqref{ineqnegativ}  we have obtained 
$$
V^{\prime}(s)\leq f(V(0)), \quad 0< s<T^*.
$$
By  a simple integration over the interval $(0,t]$ with $t<T^*$,  we  deduce that
$$
 V(t)-V(0)\leq tf(V(0)),
 $$
  for all $0\leq t<T^*$,
 which is nothing else than the first inequality of \eqref{weakbound}.
 The second inequality of \eqref{weakbound} is proved using the  inequality just above, 
 and the fact that  $f$ is non decreasing. This leads to
$$
 V^{\prime}(t)\leq f( V(0)+tf(V(0))),
 $$
 for all $0< t<T^*$.
 This completes the proof of Proposition \ref{diffineq}.
\hfill $\square$
\vskip0.3cm

\begin{rem}
In part (a)(i), under the assumption $V(0)<\lambda^*$, note that the inequality  $V^{\prime}(t)\leq  f(V(0))<0=f(\lambda^*)$
 is  a self-improvement of  $V^{\prime}(t)\leq 0$, obtained just earlier.
\end{rem}

\begin{rem}
We  have mentioned   the bound \eqref{boundTV(0)} on $T^*$ for its simplicity, and 
because this bound is used in many papers to show the finiteness of $T^*$, see for instance  \cite{CC1}.
The preceding  bound  \eqref{tcgeneral} on $T^*$ is 
a sharp   improvement of  this bound \eqref{boundTV(0)}.
\end{rem}

\begin{rem}
The argument for  the proof of (1)  in  Proposition  \ref{diffineq} is  elaborated for instance from   
the argument used by T. Coulhon to  obtain ultracontractive bounds for semigroups of operators
when the corresponding  infinitesimal  generator satisfies a generalized Nash inequality, see \cite[p.512]{coulh}.
\end{rem}

\begin{rem}
Proposition  \ref{diffineq} is sharp in the following sense. 
Let   $f$ and $\lambda^*$  be as in  Proposition \ref{diffineq}. 
For every  $V_0\in (0,\lambda^*)$, there exists a unique  function
$V:[0,T^*)\rightarrow[0,+\infty)$ such that   $T^*>0$, and $V$ 
is a non-negative continuous function, which is a  $C^1$-function  on $(0,T^*)$,
satisfying  the next conditions:
\begin{equation}\label{egalit}
 V^{\, \prime}(t)= f(V(t)), \quad 0<t<T^*, \quad  V(0)=V_0<\lambda^*,
\end{equation}
with 
\begin{equation}\label{Tmax}
T^*= \Theta_{V_0}(0^+)\in (0, +\infty].
\end{equation}
More explicitly, the unique solution $V$ is given by the following expression
\begin{equation}\label{defit}
V (t)=\Theta_{V_0}^{-1}(t),\quad 0\leq t<T^*,
\end{equation}
where 
\begin{equation}\label{deftheta}
\Theta_{V_0}(x)=\int_x^{V_0} \frac{ds}{-f(s)}, \quad  0< x\leq V_0.
\end{equation}
Hence,  inequalities \eqref{ineqdiff},  \eqref{tcgeneral} and  \eqref{thetaboundgeneral} hold true. 
In fact,  these inequalities  are equalities.  

The proof  of the assertions just above follows the same arguments as those used
in  Proposition  \ref{diffineq} in  case of equality \eqref{egalit}.
Thus , we first obtain   the relation $t=\Theta_{V_0}(V(t))$ by integration of \eqref{egalit}.
This leads to the  expression  \eqref{defit} of $V(t)$
for all $0\leq t<T^*:= \Theta_{V_0}(0^+)$, because 
the function $\Theta_{V_0}$ is  invertible.
The maximal  value of $T^*$ is defined via  the condition determined by  the range of $\Theta_{V_0}$,
$$
 {\rm Ran}\,  \Theta_{V_0}= [\Theta_{V_0}(V_0), \Theta_{V_0}(0^+) )=:[0,T^*),
$$
i.e. $T^*= \Theta_{V_0}(0^+)\in (0, +\infty]$.
\vskip0.3cm

Note that $T^*$ is finite if and only if $s\mapsto -1/{f(s)}$ is integrable at $0^+$.
Under the assumptions of  Proposition  \ref{diffineq}, this is always  the case  since $-f(0^+)\in (0,+\infty]$,
then the function $s\mapsto -1/{f(s)}$ is bounded on $[0,V_0]$.
On the other hand,  by extension of the function  $ \Theta_{V_0}$ with $V_0= \lambda^*$, 
the function $\Theta_{\lambda^*}(x)$  may be $+\infty$ at all points $x\in (0, \lambda^*)$.
Indeed,  the function $s\mapsto -1/{f(s)}$ may not be integrable at $\lambda^*$ due to the fact that $f(\lambda^*)=0$.
\end{rem}

\section{Corollaries of Theorem \ref{bupp}} \label{lescoro}
 An unexpected  consequence of Theorem \ref{bupp} is to propose a simple procedure to create blow-up solutions.
It simply consists in  multiplying any   initial data with finite variance  and positive mass by a sufficiently large constant.
Note that the center of mass and  variance are preserved under scaling of the initial datum.
See Corollary \ref{blowcreation}.

\vskip0.3cm

\subsection{Corollaries \ref{timedecay} and \ref{blowcreation}} \label{deuxcoro}
Recently, for  the case  $\alpha=0$  in equation \eqref{pks2}, the authors of \cite{Li-Wang} proved
that for any non-negative initial data  $n_0\in L^1(\R^2)$ with $M>8\pi$, 
there exist some constants $C$ and  $\lambda_0>0$ such that
$$
T^*(\lambda n_0)\leq 
\frac{C}{\lambda},
$$
 for all $\lambda\geq \lambda_0$. The second moment condition is not assume for their result.

In the case  of  finiteness of the second moment condition of  $n_0$ with mass $M>8\pi$, it can be shown that
$$
T^*(\lambda n_0)\leq 
\frac{2\pi I(0)}{(\lambda M-8\pi)M}
\leq 
\frac{C_1}{\lambda},
$$
 for all  $\lambda \geq  \lambda_0=1$ with  $C_1=\frac{2\pi I(0)}{( M-8\pi)M}$
 and
 $I(0)=\int_{\Ri^2} \vert x \vert^2 n_0(x)\, dx$.
 (See Remark 1 in \cite{Li-Wang}). 
We can refine this result  by taking into account of the centre of mass $B_0$,  which leads to
$$
T^*(\lambda n_0)\leq 
\frac{2\pi V_2(n_0)}{\lambda M-8\pi}
\leq 
\frac{C_2}{\lambda},
$$
for all $\lambda  \geq \lambda_0=1$ with   $C_2=\frac{2\pi V_2(n_0)}{M-8\pi}$.

The case $\alpha=0$ can be seen as a limit case of our estimate \eqref{Tboundzero} when $\alpha>0$ tends to $0$. 
(In that case,  the condition \eqref{eqgamma} is deleted.)
When $\alpha>0$, we  obtain a straightforward  consequence  from  Theorem \ref{bupp} 
concerning the decay of the  maximal existence  time $\lambda \mapsto T^*(\lambda n_0)$ for multiples  $\lambda\geq 1$ of a fixed   initial
data $n_0$ under finite second moment assumption.
 
 As already mentioned,  we express our results in the next corollary in terms of variance (and center of mass) 
of the initial data $n_0$  to formulate blow-up to (PKS) solutions.

 \begin{cor} \label{timedecay}
Assume that  the function $n_0$  satisfies  the same conditions as in Theorem \ref{bupp} on the mass and the second  moment.
Then, we have the next  estimate on  the maximal  existence  time $T^*(\lambda n_0)$ 
for  the solution of the (PKS) equation with initial data
$\lambda n_0$ where ${\lambda}\geq 1$,
\begin{equation}\label{Tlambda}
T^*(\lambda n_0)\leq 
\frac{1}{\lambda}
T^*_{\alpha} (n_0)
\end{equation}
with 
 $$T^*_{\alpha} (n_0):=
2\pi
\int_0^{V_2(n_0)} \frac{ds}
{
 M g_{1}\left(\sqrt{2\alpha s}\right)
-8\pi
}.
$$
\end{cor}

 {\bf Proof}. 
 Let ${\lambda}\geq 1$. Then the  mass  of $\lambda n_0$
 satisfies ${\lambda}M >8\pi$.
We also note that $V_2(\lambda n_0)=V_2(n_0)$ for all  ${\lambda}\geq 1$, i.e. the moment is invariant by considering multiples of $n_0$.
We easily check that 
 \begin{equation}\label{}
V_2(\lambda n_0)
=V_2(n_0)
<
\frac{1}{2\alpha}\left(
g_{1}^{-1}\right)^2
(8\pi / M)
\leq 
\frac{1}{2\alpha}\left(
g_{1}^{-1}\right)^2
(8\pi / \lambda M)
=
\gamma^*(\alpha, \lambda M),
\end{equation}
since  the function $g^{-1}_{1}$ is decreasing and  ${\lambda}\geq 1$.
\vskip0.3cm
Thus, we can apply Theorem \ref{bupp} to the initial data 
$n=\lambda n_0$ and obtain  the following bound on  the corresponding maximal existence time
\begin{equation} 
T^*(\lambda n_0)\leq 
2\pi
\int_0^{V_2(n_0)} \frac{ds}
{
\lambda M g_{1}\left(\sqrt{2\alpha s}\right)
-8\pi
}
\leq
\frac{2\pi} {\lambda}
\int_0^{V_2(n_0)} \frac{ds}
{
 M g_{1}\left(\sqrt{2\alpha s}\right)
-8\pi
}.
\end{equation}
 The proof is completed.
 \hfill $\square$
  \vskip0.3cm 
 
 Of course, a  similar result  can be proved using   \eqref{tcc}   and  \eqref{tks} mentioned above.
\vskip0.3cm 
 
Below, we provide an interesting variant of Corollary \ref{timedecay} by emphasizing another aspect of 
Theorem \ref{bupp},  i.e.  the  function   
  $$
  M\mapsto \gamma^*(\alpha,M):=\frac{1}{2\alpha}\left(
g_{1}^{-1}\right)^2
(8\pi / M),
$$
used to express  
the blow-up condition \eqref{eqgamma}  on  the variance $V_2(n_0)$
 is not bounded as a function of $M$. As a consequence,
we can propose a class  of initial data for which the corresponding (PKS)  solutions blow up in finite time by 
noting  that the variance  is constant for multiples of $n_0$, i.e.  it satisfies  $V_2(\lambda n_0)=V_2(n_0)$ for all $\lambda >0$.
\vskip0.3cm

Roughly speaking, the next corollary says that  for any initial data $n_1$  with finite second moment and  positive  mass,  
we have a blow-up  of all  the solutions to  the (PKS) equation with the initial data $\lambda n_1$ when  $\lambda$ is  large enough.
This simple procedure allows us to create a family of  blow-up solutions.
Such a result   does not seem to be available with  the bounds 
$\gamma^*_{cc}(\alpha,M)$ and $ \gamma^*_{ks}(\alpha,M)$ in place of $\gamma^*(\alpha,M)$,
see   \eqref{gammacc} and   \eqref{gammaks} respectively.  
Indeed, both  functions are bounded as functions of the mass $M$.
\vskip0.3cm

Here is the description of the family of initial data mentioned above.
Let $v>0$ and $b_0\in \R^2$ be fixed. 
We define  ${\mathcal E}_{v,b_0}$ as  the set of non-negative functions of $L^1(\R^2,dx)$ with common variance  $v$ 
and    center of mass $b_0$. 
(We assume that each of these functions have a  finite second moment.)
Formally,   let us define
$$
{\mathcal E}_{v,b_0}
=
\{ n\in L^1(\R^2,dx)\vert \; n\geq 0, \int_{\Ri^2} \vert x\vert^2\, n(x)\, dx<+\infty, \;  V_2(n)=v, \;B_0(n)=b_0\}.
$$
We   will   denote by $B_0(n)$  and  $V_2(n)$ the center of mass  and the variance  respectively of the function 
$n\in {\mathcal E}_{v,b_0}$ as defined in  \eqref{bary} and \eqref{varzero}.
This set is not empty, as it can be shown  by considering the function $n$ 
defined as the characteristic function of the ball of center $b_0\in \R^2$ and radius $\sqrt{2v}$.
 Note  that the mass of  any $n\in {\mathcal E}_{v,b_0}$  is non-zero since $v>0$.
 A simple  but  fundamental property of the set ${\mathcal E}_{v,b_0}$ is the following one.
For all $\lambda>0$  and all  $n\in {\mathcal E}_{v,b_0}$, 
we have  $\lambda n\in {\mathcal E}_{v,b_0}$.
In particular, the center of mass  and  the  variance of a function are invariant under  scaling $\lambda\mapsto \lambda n$,
i.e.  $B_0(\lambda n)=B_0(n)$  and  $V_2(\lambda n)= V_2(n)$ for all $\lambda>0$ and  all $n\in {\mathcal E}_{v,b_0}$. 
These  properties are crucial for the next result.
 \vskip0.3cm
 
 \begin{cor} \label{blowcreation}
 Let $m_0\in {\mathcal E}_{v,b_0}$  with  fixed variance  $v>0$ and fixed center of mass  $b_0\in \R^2$. 
 For all $\alpha>0$, let us define $ \lambda _{\alpha}(m_0)$ by
 \begin{equation}\label{lambdao}
 \lambda _{\alpha}(m_0)=
 \frac{8\pi}
 {M_0 \, g_1(\sqrt{2\alpha v})},
 \end{equation}
where $M_0=\vert\vert m_0 \vert\vert_1$.
Then, for all $ \lambda > \lambda _{\alpha}(m_0)$, the solution $(n_t)$ to the (PKS) equation with  initial data $n_0=  \lambda m_0$
 blows up in finite time.
 
 Moreover, for all $\varepsilon >0$, there exists a finite  constant  $C_{\alpha}(\varepsilon)$ such that for all
 $$
  \lambda \geq  \lambda _{\alpha}(m_0)+\varepsilon,
 $$
the maximal existence time  $T^*( {\lambda} m_0)$ of this solution satisfies
 $$
 T^*( {\lambda} m_0)
 \leq 
 \frac{
  C_{\alpha}(\varepsilon,m_0)
 }
 {\lambda}.
 $$
 In particular, we have
 $$
 \lim_{\lambda \rightarrow +\infty}  T^*( {\lambda} m_0)=0.
 $$
 More precisely, the constant $C_{\alpha} (\varepsilon,m_0)$ can be chosen  as 
\begin{equation}\label{czero}
 C_{\alpha} (\varepsilon,m_0)= 
  2\pi 
\int_0^{ v}
\frac{ds}{
 M_0 \, g_{1}\left(\sqrt{2\alpha s}\right)
-\frac{8\pi}{\lambda _{\alpha}(m_0)+\varepsilon }}.
\end{equation}
 \end{cor}
 
 {\bf Proof}. 
Fix a function $m_0\in {\mathcal E}_{v,b_0}$ with center of mass $b_0$
and variance $V_2(m_0)=v>0$. We denote by $M_0$ the mass of $m_0$, i.e. 
$M_0=\vert\vert m_0\vert\vert_1>0$.
From the assumption $ \lambda  > \lambda_{\alpha}(m_0)$
and  the expression  \eqref{lambdao}  of $ \lambda_{\alpha}(m_0)$, we  first deduce that the mass of $ \lambda  m_0$
satisfies
$\vert\vert \lambda  m_0\vert\vert_1=\lambda M_0>8\pi$,
because $0<g_1\leq 1$. The condition  \eqref{lambdao}  also implies that
$$
v< \frac{1}{2\alpha}\left(
g_{1}^{-1}\right)^2
(8\pi / \lambda M_0).
$$
(In fact, both conditions are equivalent).
Since the variance verifies  that  $V_2( \lambda  m_0)=V_2(m_0)=v$, 
then  the condition \eqref{eqgamma}  
of Theorem \ref{bupp} is satisfied with  $n_0= \lambda  m_0$, i.e.
$$
V_2( \lambda  m_0)
<
\frac{1}{2\alpha}\left(
g_{1}^{-1}\right)^2
(8\pi / \lambda M_0).
$$

Applying Theorem \ref{bupp} with initial data $n_0=\lambda  m_0$, we deduce that  the solution $(n_t)$ to the 
(PKS) equation blows up in finite time, and that the maximal existence time
$T^*(\lambda  m_0)$ of this solution satisfies  inequality  \eqref{Tboundzero}, i.e.
$$
T^*(\lambda  m_0)
\leq T^*_{\alpha} (\lambda  m_0):=
2\pi
\int_0^{V_2(\lambda  m_0)} \frac{ds}
{
\lambda M_0 g_{1}\left(\sqrt{2\alpha s}\right)
-8\pi
}
$$
$$
=
2\pi
\int_0^{v} \frac{ds}
{
\lambda M_0 g_{1}\left(\sqrt{2\alpha s}\right)
-8\pi
}
=
\frac{2\pi}
{\lambda}
\int_0^{v} \frac{ds}
{ M_0 g_{1}\left(\sqrt{2\alpha s}\right)
-\frac{8\pi}{\lambda}
}.
$$
Let $\varepsilon >0$.
Since the following  function 
$$
\lambda \mapsto \frac{1}
{ M_0 g_{1}\left(\sqrt{2\alpha s}\right)
-\frac{8\pi}{\lambda}
}
$$
 is decreasing on $(\lambda_{\alpha}(m_0),+\infty)$, then we obtain 
 for all $ \lambda \geq \lambda_{\alpha}(m_0)+\varepsilon$, the next estimate
 $$
 T^*(\lambda  m_0)
 \leq 
 \frac{2\pi}{\lambda}
\int_0^{v} \frac{ds}
{ M_0 g_{1}\left(\sqrt{2\alpha s}\right)
-\frac{8\pi}{\lambda _{\alpha}+\varepsilon}}
=:
 \frac{C_{\alpha}(\varepsilon,m_0)}
{\lambda}.
$$

This completes the proof of Corollary \ref{blowcreation}.
 \hfill $\square$
 \vskip0.3cm
 
\begin{rem}
\begin{enumerate}
\item
Corollary \ref{blowcreation} above provides  an explicit   control  by $T_{\alpha}^*( {\lambda} m_0)$  
of the value of the explosion time $T^*( {\lambda} m_0)$ for any
$m_0 \in {\mathcal E}_{v,b_0}$ with ${\lambda}$ large enough, namely ${\lambda}> \lambda_{\alpha}(m_0)$ where $\lambda_{\alpha}(m_0)$  is given by \eqref{lambdao}.
This control  is of course independent of the center of mass $b_0$, but  it depends on the variance $v>0$, 
the mass of the initial data ${\lambda} m_0$ and $\alpha>0$.
Moreover, the value of   $T^*( {\lambda} m_0)$  can be arbitrary small by choosing an  appropriate value of   $\lambda$ such that $\lambda \geq \lambda_{\alpha}(m_0)$.
In particular,  we have
$$
\lim_{ {\lambda}\rightarrow +\infty} T^*( {\lambda} m_0)=0,
$$
for each fixed $m_0 \in {\mathcal E}_{v,b_0}$ where  $v>0$ and $b_0\in \R^2$.
As a consequence, this corollary allows us to build many  solutions  to the (PKS) equation with arbitrary small explosion  time  $T^*$.
\item
Note that the initial data  ${\lambda}m_0 $ considered just above with blow-up have the same  prescribed center of mass  and  variance as 
the solution to the  (PKS) equation with   initial  data $m_0$ which globally in time exists if its mass satisfies  $\vert\vert m_0\vert\vert _1 <8\pi$ 
(no blow-up).
\item
From a heuristic point of view,  it is  not surprising that a blow-up  to the solution to the (PKS) equation  
  occurs if we increase the mass of the cells to a certain  amount,  
keeping fixed  the center of mass  and the variance of the initial data.
Corollary \ref{blowcreation}  quantifies this phenomenon.
\end{enumerate}
\end{rem}

\subsection{Corollaries  \ref{firstcoro}, \ref{secondcoro} and  \ref{explicitbup} } \label{lastcoro}

 Another interesting aspect  of Theorem \ref{bupp} is that by using several types  of estimate on   $g_1$ and  
 its inverse $g_1^{-1}$ obtained in  Section \ref{gestim}, and applying this Theorem \ref{bupp}, 
we  obtain Corollaries  \ref{firstcoro}, \ref{secondcoro} and  \ref{explicitbup}, where  explicit bounds 
on the variance $V_2(n_0)$ of the initial data $n_0$ imply the blow-up  of the solution of  the (PKS) equation.
Moreover, we   provide  explicit  upper  estimates of the maximal existence  time  $T^*$ of the solutions under the corresponding  assumptions. 

Recall that  $B_0=\int_{\Ri^2}  x\,n_t(x)\, dx/M$ 
is the center of mass  of the initial data   $n_0$, 
and $V_2(n_0)$  its variance  defined by 
$$
V_2(n_0)=\frac{1}{M} \int_{\Ri^2} \vert x-B_0\vert^2 n_0(x)\, dx
=\frac{I(0)}{M}-  \vert B_0\vert^2,
$$
where $I(0)= \int_{\Ri^2} \vert x \vert^2 n_0(x)\, dx$ is the second moment of $n_0$. 

\begin{cor}\label{firstcoro}
Assume that $(n_t)_{0<t<T^*}$ is a solution of  (PKS) equation \eqref{ks} for a fixed $\alpha>0$.
We denote by  $T^*$ the maximal  existence time  
of the solution  with initial condition $n(0,x)=n_0(x)\geq 0$, where $n_0\in L^1(\R^2)$.
We assume that $M=\int _{\Ri^2} n_0(x)\, dx>8\pi$, and the second moment $I(0)$ of $n_0$ is finite.
If the variance satisfies the following bound
\begin{equation}\label{eqgamma**}
V_2(n_0)
<
\gamma^{**}(\alpha, M)
:=
\frac{1}{2\alpha}
\ln^2\left( \frac{M}{8\pi}
\right),
\end{equation}
then the solution $(n_t)$ of the (PKS) equation blows up in finite time, i.e. $T^*<+\infty$,
and   we  have the following bounds on $T^*$,
 \begin{equation}\label{Tboundzero**}
T^*\leq 
T^{**}_c(\alpha, n_0)
:
=
2\pi
\int_0^{V_2(n_0)} \frac{ds}
{
M  e^{-\sqrt{(2\alpha \, s)}}
-8\pi
}
\leq
L(\alpha,n_0)
\end{equation}
where
 \begin{equation}\label{Lnzero}
L(\alpha,n_0)
:=
\frac{2\pi\, V_2(n_0)}
{
M  e^{-\sqrt{2\alpha \, V_2(n_0)}}
-8\pi
}.
\end{equation}
\end{cor}

Note that the bound $\gamma^{**}(\alpha, M)$ in \eqref{eqgamma**}
depends only on  the mass of the initial condition $n_0$, but the bound $T^{**}_c(\alpha, n_0)$ on $T^*$ depends  also on the  variance $V_2(n_0)$
of the initial data $n_0$ for  fixed $\alpha>0$.

From (4) of Proposition \ref{propgun}, we have  that  
$g_1^{-1}(\frac{8\pi}{M})\sim \ln(\frac{M}{8\pi} )$ as $M$ tends to $+\infty$.
Thus, the assumption \eqref{eqgamma**} on the variance   is asymptotically sharp with respect to the assumption
\eqref{eqgamma} of Theorem \ref{bupp} when $M$ tends to infinity.
\vskip0.3cm

{\bf Proof}. 
By inequality  \eqref{uniflowerginv} of Lemma \ref{lowerg1}, we have for all $\rho \in (0,1]$,
$$
\ln\left(\frac{1}{\rho}\right)\leq 
g_1^{-1}(\rho).
$$
Taking $\rho=\frac{8\pi}{M}$ when $M>8\pi$, we deduce that
$$
\frac{1}{2\alpha}
\ln^2\left( \frac{M}{8\pi}
\right)
\leq
\frac{1}{2\alpha}\left(
g_{1}^{-1}\right)^2
(8\pi /M),
$$
for any $\alpha >0$.
Hence, the assumption  \eqref{eqgamma**} implies the assumption \eqref{eqgamma}.
By applying  \eqref{Tboundzero} of Theorem  \ref{bupp},
we immediately deduce that
 $$
T^*\leq T^*_c(\alpha, n_0):=
2\pi
\int_0^{V_2(n_0)} \frac{ds}
{
M g_{1}\left(\sqrt{2\alpha s}\right)
-8\pi
}
.
$$
In particular, we have $T^*<+\infty$.
Now, using  the exponential-type   lower bound  \eqref{uniflowerg}  on $g_1$ obtained in  Lemma \ref{lowerg1},
i.e. $e^{-r}\leq g_1(r)$ for any $r\geq 0$, 
 and using the fact that 
 $$
 s\mapsto 
\frac{1}
{
M  e^{-\sqrt{(2\alpha \, s)}}
-8\pi
}
$$
 is increasing,  we finally obtain successively 
$$
2\pi
\int_0^{V_2(n_0)} \frac{ds}
{
M g_{1}\left(\sqrt{2\alpha s}\right)
-8\pi
}
\leq
2\pi
\int_0^{V_2(n_0)} \frac{ds}
{
M  e^{-\sqrt{(2\alpha \, s)}}
-8\pi
}
$$
$$
\leq
\frac{2\pi\, V_2(n_0)}
{
M  e^{-\sqrt{2\alpha \, V_2(n_0)}}
-8\pi
}:=L(\alpha,n_0).
$$
The proof of the corollary is completed.
\hfill $\square$
\vskip0.3cm

In the next result, we obtain   an expression  in terms of series of  the integral $T^{**}_c(\alpha, n_0)$ 
appearing in \eqref{Tboundzero**} of  Corollary \ref{firstcoro}.
 We also find another expression in terms of well-known polylogarithms defined  for $s\in \C$ by the series
$$
{\rm Li}_s(x)=\sum_{n\geq 1}\frac{x^n}{n^s},\; \vert x\vert<1.
$$

Furthermore, we provide below another  upper estimate  $K(\alpha, n_0)$  
of this integral to be compared with $L(\alpha, n_0)$ of Corollary \ref{firstcoro}.

\begin{cor}\label{secondcoro}
Under the same assumptions  as in  Corollary \ref{firstcoro}.
Let  $T^{**}_c(\alpha, n_0)$ be the critical time given by the expression
\begin{equation}\label{critic**}
T^{**}_c:=T^{**}_c(\alpha, n_0)
=
2\pi
\int_0^{V_2(n_0)} \frac{ds}
{
M  e^{-\sqrt{(2\alpha \, s)}}
-8\pi
}, 
\end{equation}
which satisfies $T^*\leq T^{**}_c(\alpha, n_0)$.
Then the following estimates hold true,
\begin{enumerate}
\item
\begin{equation}\label{Tseries}  
T^{**}_c=
\frac{2\pi}{\alpha}
\sum_{n=0}^{\infty}
\frac{1}{(n+1)^2}
\frac{(8\pi)^n}
{M^{n+1}}
\left[1+
\left( (n+1)\sqrt{ 2\alpha V_2(n_0)}
-1\right) e^{
 (n+1)\sqrt{ 2\alpha V_2(n_0)}
}
\right].
\end{equation}
\vskip0.3cm

\item 
\begin{equation}
T^{**}_c
=
\ell(\alpha)
\sqrt{V_2(n_0)}
\ln\left(
\frac{M }
{M-8\pi e^{\sqrt{2\alpha \, V_2(n_0)}}}
\right)+U(n_0)
 \end{equation}
with
\begin{equation}\label{upperlis} 
U(n_0)=
-\frac{1}{4\alpha}{\rm Li}_2\left(
\frac{8\pi  }
{M}
e^{\sqrt{2\alpha \, V_2(n_0)}}
\right)
+\frac{1}{4\alpha}
{\rm Li}_2\left(
\frac{8\pi  }
{M}\right) 
< 0,
 \end{equation}
where $\ell(\alpha)=\frac{1}{2 \sqrt{2\alpha}}$, $\alpha >0$.
\vskip0.3cm

\item 
\begin{equation}\label{upperTstar} 
T^{**}_c
\leq
K(\alpha, n_0):=
\ell(\alpha)
\sqrt{V_2(n_0)}
\ln\left(
\frac{M-8\pi }
{M-8\pi e^{\sqrt{2\alpha \, V_2(n_0)}}}
\right).
\end{equation}
\vskip0.3cm

\item (First comparison between $K(\alpha, n_0)$ and $L(\alpha, n_0)$).

Let $Y_0$ be the unique value such that 
$Y_0>1$  and  $(\frac{1}{2}Y_0-1) e^{Y_0}+1=0$
($Y_0\sim 1,594$).  Assume that 
\begin{equation}\label{yzerocompar}
Y_0< \sqrt{2\alpha \, V_2(n_0)} <\ln\left(\frac{M}{8\pi}\right),
\end{equation}
then we have 
\begin{equation}\label{compare lk}
\ell(\alpha)
\sqrt{V_2(n_0)}
\ln\left(
\frac{M-8\pi }
{M-8\pi e^{\sqrt{2\alpha \, V_2(n_0)}}}
\right)
<
\frac{2\pi\, V_2(n_0)}
{
M  e^{-\sqrt{2\alpha \, V_2(n_0)}}
-8\pi
},
\end{equation}
i.e.
$K(\alpha, n_0) <L(\alpha,n_0)$,
where $L(\alpha,n_0)$ and $K(\alpha, n_0)$ 
are given respectively by  \eqref{Tboundzero**} and \eqref{upperTstar}.
As a consequence, the estimate  \eqref{upperTstar} of $T^*$ improves   the estimate    \eqref{Tboundzero**}
of $T^*$  as soon as \eqref{yzerocompar} is verified.
\vskip0.3cm

\item (Second comparison between $K(\alpha, n_0)$ and $L(\alpha, n_0)$).

Assume that $n_0$ is  such that   $M>8\pi$ and
$V_2(n_0)<\frac{1}{2\alpha} \ln^2( M/8\pi)$.
Let $Y_2=Y_2(M)$ be the unique solution of the next equation,
$$
Y-1+\frac{8\pi}{M} e^Y=0,
\quad Y\in  (0,  \ln( M/8\pi)),
$$
 and let 
 $Y_1=Y_1(M)\in(Y_2,  \ln( M/8\pi))$
be the unique solution of
$$
\frac{4\pi Y e^{Y}}
{M-8\pi  e^{Y}}
=
\ln\left(
\frac{M-8\pi}
{M-8\pi e^{Y}}\right),
\quad Y\in  (0,  \ln( M/8\pi)).
$$
Then, 
 $$
 K(\alpha, n_0) \leq L(\alpha, n_0)
 $$
 if and only if if $Y_1\leq  \sqrt{2\alpha V_2(n_0)} <\ln(M/8\pi)$.

 (The case of equality is obtained if and only if  $Y_1=\sqrt{2\alpha V_2(n_0)}$).
 \end{enumerate}
\end{cor}
 
Note that the bounds   $L(\alpha, n_0)$ and $K(\alpha, n_0)$  defined  as in     
Corollary  \ref{firstcoro}  and Corollary  \ref{secondcoro} respectively,
cannot be bounded  by a function depending only  on  the variance of $n_0$.
Indeed, the quantity ${M-8\pi e^{\sqrt{2\alpha \, V_2(n_0)}}}$ can be arbitrary small,
for instance by  fixing  the variance and letting $M\rightarrow 8\pi e^{\sqrt{2\alpha \, V_2(n_0)}}$
 using invariance  by scaling of the variance  of $n_0$.
 \vskip0.3cm
 
Statement (4) of Corollary \ref{secondcoro} characterizes  the fact that  $K(\alpha, n_0)$ is a better bound of $T^*$  than
 $L(\alpha, n_0)$ depending on whether or not  the variance satisfies the condition $Y_1(M)< \sqrt{2\alpha V_2(n_0)}$.
 Note that $Y_1(M)$ and $Y_2(M)$ depend only on  the mass $M$, and not on the variance of $n_0$.
It is easy to see that we have $0< Y_2<1$ because 
 $$
\frac{8\pi}{M} e^{Y_2}=1-Y_2>0.
$$
We can improve this estimate of $Y_2$ by showing the next lower and upper bounds:
$$
B_1(M):=\ln\left(\frac{2M}{M+8\pi}\right)
 < Y_2(M),
 $$
 and 
 $$
Y_2(M) <
B_2(M):=
\frac{1}{2}\left[
\sqrt{
\left(\frac{M}{4\pi}+2\right)^2+\left(\frac{M}{\pi} -8\right)
}
-\left(\frac{M}{4\pi}+2\right)
\right].
$$
The lower bound is proved by studying the inequality 
$$
(1+\frac{8\pi}{M}) e^{Y_2} >1+{Y_2}+\frac{8\pi}{M} e^{Y_2}=2,
$$
and the upper bound by studying the quadratic  inequality 
$$
1-{Y_2}= \frac{8\pi}{M} e^{Y_2}> \frac{8\pi}{M} \left(1+{Y_2}+\frac{{Y_2}^2}{2}\right).
$$
Note that $M\in (8\pi,+\infty)\mapsto Y_2(M)$ is well defined since the solution $Y_2$ is unique for each $M>8\pi$, 
and it is a smooth function  (strictly increasing in $M$)  by using  implicit function theorem.
\vskip0.3cm
 
{\em Two mass examples}.
\begin{enumerate}
\item
Let $M=16\pi$. Then, we have $B_1(M)\sim 0,288$ and $B_2(M)=0.316$.
\item
Let $M=24\pi$. Then,  we have $B_1(M)=0.405$ and $B_2(M)=0.472$.
\end{enumerate}
 For any fixed  $M>8\pi$, numerical solvers can be used to provide accurate approximations of   $Y_i$'s values.
 For the first case, i.e. $8\pi/M=1/2$,  we have $Y_1\sim 0.461$ and  $Y_2\sim  0.315$.
For the second case, i.e. $8\pi/M=1/3$, we have $Y_1\sim 0.693 $ and  $Y_2\sim 0.468$.
Note that the value $B_2(M)$ is closed to the true solution $Y_2(M)$ in  both examples.
\vskip0.3cm

{\bf Proof}.
(1) By definition, we have
$$
T^{**}_c(\alpha, n_0)
:
=
2\pi
\int_0^{V_2(n_0)} \frac{ds}
{
M  e^{-\sqrt{2\alpha \, s}}
-8\pi
}.
$$ 
Let $X= V_2(n_0)$. Then we can write
$$
T^{**}_c(\alpha, n_0)
=
2\pi
\int_0^X \frac{e^{\sqrt{2\alpha \, s}}}
{
M 
-8\pi
 e^{\sqrt{2\alpha \, s}}
}
\, ds
=
\frac{2\pi}{M}
\int_0^X \frac{e^{\sqrt{2\alpha \, s}} }
{1
-\frac{8\pi}{M}
 e^{\sqrt{2\alpha \, s}}
}
\, ds
$$ 
$$
=
\frac{2\pi}{M}
\int_0^X e^{\sqrt{2\alpha \, s}} 
\sum_{n=0}^{\infty}
\left( \frac{8\pi}{M}\right)^n 
e^{n\sqrt{2\alpha \, s}}
=
\frac{2\pi}{M}
\sum_{n=0}^{\infty}
\left( \frac{8\pi}{M}\right)^n 
\int_0^X 
e^{(n+1)\sqrt{2\alpha \, s}}
ds.
$$
Here,  we have used Fubini's Theorem in the last equality, and for the series,  the fact  that
$$
0< \frac{8\pi}{M}
 e^{\sqrt{2\alpha \, s}}
 \leq 
 \frac{8\pi}{M}
  e^{\sqrt{2\alpha \, V_2(n_0)}}
 <1,
 $$ 
 for all $s\in [0, V_2(n_0)]$. The last  inequality is exactly the assumption above on the variance 
 $V_2(n) <\frac{1}{2\alpha} \ln^2(M/8\pi)$.
 \vskip0.3cm

By the change of variables $t= (n+1)\sqrt{2\alpha \, s}$, we obtain 
$$
T^{**}_c(\alpha,n_0)
=
\frac{2\pi}{\alpha M}
\sum_{n=0}^{\infty}
\frac{1}{(n+1)^2}
\left( \frac{8\pi}{M}\right)^n 
\int_0^{y_n}
te^{t}\,
dt,
$$
with $y_n=(n+1)\sqrt{2\alpha \, X}$. The last integral can be computed explicitly by integration by parts. 
We finally get
$$
T^{**}_c(\alpha,n_0)
=
\frac{2\pi}{\alpha M}
\sum_{n=0}^{\infty}
\frac{1}{(n+1)^2}
\left( \frac{8\pi}{M}\right)^n 
\left[
(y_n-1)e^{y_n}+1
\right].
$$
So, the formula \eqref{Tseries} is proved.
\vskip0.3cm

(2)
 Let   $X= 2\alpha V_2(n_0)$ and  $\ell(\alpha)=\frac{1}{2 \sqrt{2\alpha}}$ for any fixed $\alpha >0$.
By the change of variables  $v=2\alpha s$,  we can write
$$
T^{**}_c(\alpha, n_0)
=
2\pi
\int_0^{V_2(n_0)} 
\frac{1}
{Me^{-\sqrt{2\alpha \, s}} -8\pi }
\, ds
=
\frac{1}{8\alpha}
\int_0^{X} \frac{1}{a e^{-\sqrt{v}}-1} \, dv
=:
\frac{1}{8\alpha} J(X),
$$
where $a =\frac{M}{8\pi}>1$.
By the new change of variables $y=\sqrt{v}$ in $J(X)$, we get
$$
J(X)=
\int_0^{\sqrt{X}}
 \frac{2y}
{a e^{-y}-1}
\, dy
=
\int_0^{\sqrt{X} } 2y \left( \frac{ e^{y}}
{
a -e^{y}
}\right)
\, dy.
$$
By integration by parts, we obtain
$$
J(X)=
\left[ {-2y} \ln(a-e^y) \right]^{\sqrt{X}}_0
+ 2
\int_0^{\sqrt{X}}
\ln( a-e^y)
\, dy
=
-2\sqrt{X}
 \ln(a-e^{\sqrt{X}})
 $$
 $$
+2\sqrt{X}\ln a
+2\int_0^{\sqrt{X}}
\ln\left(1-\frac{e^y}{a}\right)
\, dy.
$$

Under the  following conditions $\frac{e^{\sqrt{X}}}{a} <1$  on $X$, we have  the next well-known expression in terms of  power series,
$$
\ln\left(1-\frac{e^y}{a}\right)=- \sum_{n\geq 1}
\frac{
(\frac{e^{ny}}{a^n})
}{n},
$$
for all  $y\in [0,\sqrt{X}]$.
Hence, we get
$$
\int_0^{\sqrt{X}}
\ln\left(1-\frac{e^y}{a}\right)
\, dy
=
- \sum_{n\geq 1}
\frac{1} {n^2}
\left(\frac{e^{n\sqrt{X}}-1}
{a^n}\right)
=
-{\rm Li}_2\left(\frac{e^{\sqrt{X}}}{a} \right)+{\rm Li}_2\left(\frac{1}{a}\right)
$$
with $a>1$.
Thus, we conclude that
$$
J(X)=
-2\sqrt{X}
 \ln\left(1-\frac{e^{\sqrt{X}}}{a}\right)
-2{\rm Li}_2\left(\frac{e^{\sqrt{X}}}{a}
\right)+ 2{\rm Li}_2\left(\frac{1}{a}\right).
$$
Replacing $a$ and $X$ by the value above in the equality 
$
T^{**}_c(\alpha, n_0)
=
\frac{1}{8\alpha} J(X),
$
we obtain formula \eqref{upperlis} with $U(n_0)<0$ because $x\mapsto {\rm Li_2}(x)$ is strictly increasing.

\vskip0.3cm
(3)
 Let   $X= V_2(n_0)$ and  $\ell(\alpha)=\frac{1}{2 \sqrt{2\alpha}}$ for any fixed $\alpha >0$.
We can write
$$
T^{**}_c(\alpha,n_0)
=
2\pi
\int_0^X \frac{e^{\sqrt{2\alpha \, s}}}
{
M 
-8\pi
 e^{\sqrt{2\alpha \, s}}
}
\, ds
= -
 \ell(\alpha)
\int_0^X 
\sqrt{s}\;
\frac{\frac{d}{ds} (M-8\pi e^{\sqrt{2\alpha \, s}} )}
{
M 
-8\pi
 e^{\sqrt{2\alpha \, s}}
}
\, ds
$$
$$
= -  \ell(\alpha)
\int_0^X 
\sqrt{s}\;
\frac{d}{ds} \ln (M-8\pi e^{\sqrt{2\alpha \, s}} )
\, ds,
$$
because $M-8\pi e^{\sqrt{2\alpha \, s}}>0$ on the integration  interval  $[0,X]$.
 Now, by integration by parts we get
$$
T^{**}_c(\alpha, n_0)
=
 -\ell(\alpha)
\sqrt{X}\;
 \ln \left(
 M-8\pi e^{\sqrt{2\alpha \, X}} \right)
+
 \ell(\alpha)
\int_0^X 
 \ln \left(M-8\pi e^{\sqrt{2\alpha \, s}} \right)
 d(s^{1/2})
$$
$$
\leq
 -\ell(\alpha)
\sqrt{X}\;
 \ln \left(
 M-8\pi e^{\sqrt{2\alpha \, X}} \right)
+
 \ell(\alpha)
  \ln (M-8\pi e^{0} )
 \int_0^X 
 d(s^{1/2}).
$$
$$
=
 -\ell(\alpha)
\sqrt{X}\;
 \ln \left(
 M-8\pi e^{\sqrt{2\alpha \, X}} \right)
+
 \ell(\alpha)
 \sqrt{X}\;  \ln (M-8\pi)
.
$$
This proves inequality \eqref{upperTstar}.
\vskip0.3cm

(4)
Here, we compare both estimates  $K(\alpha, n_0)$ and $L(\alpha,n_0)$ of $T^{**}_c$
appearing in   inequalities \eqref{Tboundzero**} and \eqref{upperTstar}  as follows.
Again, let  $X= V_2(n_0)$.
Under the conditions $Y_0< Z=\sqrt{2\alpha \, V_2(n_0)} <\ln(\frac{M}{8\pi})$, where $Y_0$ satisfies  $Y_0>1$, 
and  $\frac{1}{2}Y_0e^{Y_0}-e^{Y_0}+1=0$,
 then  we have 
$$
K(\alpha, n_0)=
\ell(\alpha)
\sqrt{X}
\ln\left(
\frac{M-8\pi }
{M-8\pi e^{\sqrt{2\alpha X\,}}}
\right)
=
\frac{1}{2 \sqrt{2\alpha}}
\sqrt{X}
\ln\left(
 1+ 
\frac{8\pi(e^Z-1)}
{M-8\pi e^Z}
\right)
$$
$$
\leq
\frac{1}{2 \sqrt{2\alpha}}
\sqrt{X}
\left(
\frac{8\pi(e^Z-1)}
{M-8\pi e^Z}
\right)
<
\frac{2\pi\, X e^Z}
{M 
-8\pi e^Z}
=
\frac{2\pi\, X}
{
M  e^{-\sqrt{2\alpha \, X}}
-8\pi
}.
$$
The last inequality is due to the fact that $e^Z-1 < \frac{1}{2}Ze^Z$ for any $Z>Y_0$.
\vskip0.3cm

Hence,  inequality  \eqref{upperTstar} improves the bound on  $T^{**}_c$ given by \eqref{Tboundzero**} 
when $V_2(n_0)$ satisfies  the conditions  $Y_0< \sqrt{2\alpha \, V_2(n_0)} <\ln(\frac{M}{8\pi})$.

\vskip0.3cm

(4) The inequality $K(\alpha, n_0)<L(\alpha, n_0)$ is equivalent to $H(Y)>0$, where
$Y=\sqrt{2\alpha \, V_2(n_0)}$ and
$$
H(Y)=\frac{4\pi Y e^Y}{M-8\pi e^Y}
-
\ln\left(\frac{M-8\pi }{M-8\pi e^Y}\right), \quad Y\in (0, \ln(M/8\pi)).
$$
We have $H(0^+)=0$ and $H(\ln(M/8\pi)^+)=+\infty$.
The derivative of $H$ is given by
$$
H^{\prime}(Y)= \frac{4\pi M e^Y}{(M-8\pi e^Y)^2} R(Y)
,
$$
where $R(Y)=\left(Y-1+\frac{8\pi}{M} e^Y\right)$. 
The function $R$ defined on the interval $(0, \ln(M/8\pi))$ is a  strictly increasing continuous function
with $R(0^+)=\frac{8\pi}{M} -1<0$, and $R( \ln(M/8\pi))= \ln(M/8\pi)>0$ ($Y=\ln(M/8\pi)$ is a fixed point for $R$).
Then, there exits a unique $Y_2=Y_2(M)\in (0, \ln(M/8\pi))$ 
such that $R(Y_2)=0$. Moreover,  if $Y\in (0, Y_2)$  then  we have 
$R(Y) <0$,   and  if $Y\in (Y_2,  \ln(M/8\pi))$ then we have
$R(Y) >0$.
\vskip0.3cm

Thus, the function $H$ is strictly decreasing on the interval $(0, Y_2)$,
and  strictly increasing on the interval  $(Y_2,  \ln(M/8\pi)))$.
Since $H$ is continuous, there exists a unique solution 
$Y_1=Y_1(M) \in  (Y_2,  \ln(M/8\pi))$  of the equation $H(Y)=0$. Moreover,
  if $Y\in (0, Y_1)$ then  
 $H(Y)<0$,  and if $Y\in (Y_1,  \ln(M/8\pi))$ then  $H(Y)>0$. We apply these results with  $Y=\sqrt{2\alpha V_2(n_0)}$ and conclude Statement (4).
\vskip0.3cm

This finishes the proof of the corollary.
\hfill
$\square$
\vskip0.3cm

Another   condition of blow-up can be described as follows.

\begin{cor}\label{explicitbup}
Let $c_{\varepsilon} =
(1- {\varepsilon} ) \sqrt{\frac{\pi}{2}}$ for any fixed ${\varepsilon} \in (0, 1)$.
Then, there exists $M_{\varepsilon}>8\pi$ such that, for any  initial data $n_0$ of the  (PKS) equation
with $M:=\vert\vert n_0\vert\vert_1 >M_{\varepsilon}$ and  satisfying 
\begin{equation}\label{eqgammastarstar}
V_2(n_0)
<
\Gamma^{*}_{\varepsilon}(\alpha, M),
\end{equation}
where
\begin{equation}\label{GAMMA}
\Gamma^{*}_{\varepsilon}(\alpha, M)=
\frac{1}{2\alpha}
\left[
\ln\left(\frac{c_{\varepsilon} M}{8\pi}
\sqrt{
\ln\left(
\frac{c_{\varepsilon} M}{8\pi}\right)}
\right)
\right]^2
.
\end{equation}
Then, the solution $(n_t)$ of the (PKS) equation  with initial data $n_0$ blows up in finite time,
i.e. $T^*<+\infty$, and  the   bounds   \eqref{Tboundzero} and  \eqref{Tbound} hold true for $T^*$.
\end{cor}
 
{\bf Proof}.
 Let ${\varepsilon} \in (0,1)$ and  $c_{\varepsilon} =(1- {\varepsilon} ) \sqrt{\frac{\pi}{2}}$.
By applying the lower bound of  inequality \eqref{g1inverse}  of Proposition \ref{propgun} below, 
 there exists $M_{\varepsilon}>8\pi$ such that, for  all $M>M_{\varepsilon}$,
$$
\Gamma^{*}_{\varepsilon}(\alpha, M)=
\frac{1}{2\alpha}
\left[
\ln\left(\frac{c_{\varepsilon} M}{8\pi}
\sqrt{
\ln\left(
\frac{c_{\varepsilon} M}{8\pi}\right)}
\right)
\right]^2
\leq \frac{1}{2\alpha}\left(
g_{1}^{-1}\right)^2
(8\pi /M).
$$
Hence,  the assumption
$
V_2(n_0)
<
\Gamma^{*}_{\varepsilon}(\alpha, M)
$
implies  \eqref{eqgamma}. The conclusion follows from Theorem \ref{bupp}.
 \hfill $\square$

\subsection{Estimates of $g_1$ and $g_1^{-1}$ } \label{gestim}
In this section, we prove several  estimates for $g_1$ and its inverse 
$g_1^{-1}$ useful for applications. 
We first start with a simple  lemma providing  an explicit uniform lower bound for  the function $g_1$.
We then  deduce the corresponding one for  the inverse function $g_1^{-1}$. 
 Details of proofs  are given, except for the well-known asymptotic of the function $g_1$
(references are provided).
\begin{lem}\label{lowerg1}
Let  $g_1$ be defined  by
\begin{equation}\label{g1def}
 g_{1}(r)=\int_0^{+\infty}e^{-\frac{r^2}{4s}} e^{-s}\, ds,
\quad r\geq 0.
\end{equation}
(See \eqref{galphafn} with $\alpha=1$). 
\begin{enumerate}
\item
For all $r\geq 0$, we have
\begin{equation}\label{uniflowerg}
e^{-r}\leq g_1(r).
\end{equation}
\item
For all $\rho \in (0,1]$, we have
\begin{equation}\label{uniflowerginv}
\ln\left(\frac{1}{\rho}\right)\leq 
g_1^{-1}(\rho).
\end{equation}
\end{enumerate}
\end{lem}

{\bf Proof}.
(1) Let  $r>0$ be fixed and  let $\lambda>0$ be chosen later. We have the following  simple inequalities
$$
g_{1}(r)
\geq 
\int_{\lambda r}^{+\infty}
e^{-\frac{r^2}{4s}} e^{-s}\, ds
\geq
e^{-\frac{r}{4\lambda}}
\int_{\lambda r}^{+\infty}
e^{-s}\, ds
=
e^{-\frac{r}{4\lambda }} 
 e^{-\lambda r}
 =
e^{-r(\frac{1}{4\lambda }+\lambda)}.
$$
This implies \eqref{uniflowerg} by taking  $\lambda=1/2$ (optimizing the right-hand side of this inequality over  $\lambda>0$).
\vskip0.3cm

(2) Let $\omega(r)=e^{-r}$ for all $r\geq 0$. 
The functions  $\omega$ and $g_1$ are both   continuous  and strictly decreasing on $[0,+\infty[$ with the same range
$(0,1]$. They also satisfies  $\omega(0)= g_1(0)=1$
and $\omega(+\infty)= g_1(+\infty))=0$. The function $\omega$ is explicitly invertible with
$\omega^{-1}(\rho)=\ln(1/\rho)$  for all $\rho\in (0,1]$.
Since the inequality $\omega(r)\leq g_1(r)$ holds for all $r\geq0$ by  \eqref{uniflowerg},
then this immediately  implies \eqref{uniflowerginv} for the inverse functions.
 The proof of the  lemma is completed.
 
 \hfill $\square$

The next proposition  provides a  precise  behavior  of  the functions $g_1$ and $g_1^{-1}$,
not only  asymptotically but also by pointwise estimates.
In fact, it is more desirable  to have  pointwise bounds  than asymptotics
in order to propose   explicit upper bound assumptions on the variance 
$V_2(n_0)$,  and blow-up bounds  for  (PKS) solutions. 
These estimates are used in some  of the corollaries of  Theorem \ref{bupp} in Section \ref{lescoro}.

\begin{pro}\label{propgun}  
Let $g_{1}$ be defined as in \eqref{g1def}. We then  have the following  statements.
\begin{enumerate}
\item  The following asymptotic holds:
$g_1(r)\sim \sqrt{\frac{\pi r}{2}}e^{-r}$ as $r$ tends to $+\infty$.
\vskip0.3cm
\item 
Let  $c>0$ and 
 $v_c(r)=c\sqrt{r} e^{-r}, r>0$. Then $v_c$ is a decreasing bijection from $(1/2,+\infty)$ onto $(0,v_c(1/2))$. 
 For any $\rho\in (0,v_c(1))$, 
we have 
\begin{equation}\label{minmaxvc}
 \ln\left( \frac{c}{\rho}\sqrt{ \ln ( \frac{c}{\rho})}\right)
 \leq
 v_c^{-1}(\rho)
\leq
\ln\left(\frac{c}{\rho} \sqrt{ \ln(\frac{c}{\rho})}\right)
+ \ln\left( \sqrt{\frac{2e}{2e-1}}\right).
\end{equation}
In particular, the asymptotics 
$v_c^{-1}(\rho)\sim \ln\left( \frac{c}{\rho}\sqrt{ \ln ( \frac{c}{\rho})}\right)
\sim \ln\left( \frac{1}{\rho}\right)$ hold true as $\rho$ goes to $0$.
\vskip0.3cm
\item
For all $\varepsilon \in (0,1)$, there exists $\rho_\varepsilon \in (0,1)$ such that, for all 
$\rho\in (0,\rho_\varepsilon)$, we have
\begin{equation}\label{g1inverse}
0<
\ln\left(\frac{c_{-}}{\rho} \sqrt{ \ln(\frac{c_{-}}{\rho})}\right)
\leq 
g_1^{-1}(\rho)
\leq
\ln\left(\frac{c_{+}}{\rho} \sqrt{ \ln(\frac{c_{+}}{\rho})}\right)
+ \ln\left( \sqrt{\frac{2e}{2e-1}}\right),
\end{equation}
where $c_{\pm}=(1\pm \varepsilon)\sqrt{\frac{\pi}{2}}$.
\vskip0.3cm
\item
The following asymptotic holds true
$g_1^{-1}(\rho)\sim \ln(\frac{1}{\rho})$ as $\rho$ tends to $0$.
\vskip0.3cm
\item
We have 
$v^{-1}(g_1(r))\sim r$ as $r$ goes to $+\infty$, and   $g_1(v^{-1}(\rho))\sim \rho $ as $\rho$ goes to $0$, 
where  $v(r)=\sqrt{\frac{\pi r}{2}}e^{-r}$ and $v^{-1}$ its inverse.  
The result is also valid when we replace $v^{-1}$ by 
$w(\rho)=\ln\left(\frac{c}{\rho} \sqrt{ \ln(\frac{c}{\rho})}\right)$ with $c=\sqrt{\frac{\pi}{2}}$ and
$\rho\in (0,  c/e)$.
 But the result is no longer true
if  we replace $v^{-1}$ by $h(\rho)=\log(1/\rho)$, because 
$h(g_1(r))\sim r$ as $r$ goes to $+\infty$,  but  $g_1(h(\rho))\nsim \rho $ as $\rho$ goes to $0$.
 \end{enumerate}
\end{pro}

We provide proofs except for claim (1) for which several references are available. 
Statement (5) is intended to justify the sharpness of pointwise bounds on $g_1^{-1}$ given in (3),
and  also to understand the behavior of $g_1^{-1}$ not only 
by  its asymptotic behavior but by upper and lower estimates. 

These estimates  are  used in Section \ref{lescoro} to propose some   explicit upper bound conditions on 
the variance $V_2(n_0)$ instead of  the condition \eqref{eqgamma}, which is  expressed in terms of the function $g_1^{-1}$
 (sharp for large masses).
We also derive  bounds on the maximal  existence time  $T^*$   to the blow-up of the (PKS) solutions.
\vskip0.3cm 

{\bf Proof.}
(1) The function $g_1$ is related to the modified Bessel function of the third kind by the formula
$g_1(r)=rK_{-1}(r)$, see  \cite[8.432(6).p.916]{GR}.  
The asymptotic form $g_1(r)\sim \sqrt{\frac{\pi r}{2}} e^{-r}$  as   $r$ tends to $\infty$ is well-known.
For instance, see \cite[p. 415 (3,5)]{A-S}, or 
  \cite[8.451(6).p.920]{GR}. 
We can also see  \cite[p.4874]{Ryznar} and  for a sketch of the proof  
we can look at  \cite[p.296. Chap.56]{Don}.
From this asymptotic, we   will  derive  in (4) the asymptotic of the inverse function: 
$g_1^{-1}(\rho)\sim \ln(\frac{1}{\rho})$ as $\rho$ tends to $0$.
\vskip0.3cm

(2) Let  $v:= v_c$ with $v_c$ defined as in  (2) of Proposition \ref{propgun}.
It is easy to see that the function 
$$
v:(1/2,+\infty)\rightarrow (0, v(1/2))= (0, {c}(2e)^{-1/2})
$$
is a continuous strictly decreasing bijection.
Let $\rho\in (0,\rho_0)$ with $\rho_0=\frac{c}{e}=v(1)<v(1/2)$. 
Let also $r=v^{-1}(\rho)$, hence $r>1$ and  $\rho=v(r)=c\sqrt{r} e^{-r}$. 
From this last equation, we deduce that
\begin{equation}\label{rho1}
r=\ln(\frac{c}{\rho})+\frac{1}{2}\ln (r).
\end{equation}
Hence, we have
$$
1<\ln(\frac{c}{\rho})=r(1-\frac{1}{2r}\ln (r)),
$$
and
\begin{equation}\label{rho2}
\ln r=\ln \ln(\frac{c}{\rho})- \ln (1-\frac{1}{2r}\ln (r)).
\end{equation}
(Note that $1-\frac{1}{2r}\ln (r)>0$ for all $r>0$, and  $\frac{c}{\rho}>1$ since $0<\rho<\rho_0=\frac{c}{e}<c$).
Combining both equations \eqref{rho1} and  \eqref{rho2} above, 
we thus obtain 
$$
r=\ln(\frac{c}{\rho})+\frac{1}{2}\ln \ln(\frac{c}{\rho})-  \frac{1}{2}\ln \left(1-\frac{1}{2r}\ln (r)\right),
$$
i.e.
\begin{equation}\label{vcequat}
v^{-1}(\rho) =r=\ln\left(\frac{c}{\rho} \sqrt{ \ln(\frac{c}{\rho})}\right)-  \frac{1}{2}\ln \left(1-\frac{1}{2r}\ln (r)\right).
\end{equation}
For  all  $r>1$,  we have  $0<\frac{\ln r}{r}\leq e^{-1}$, then this implies that
\begin{equation}\label{eee}
0<-  \frac{1}{2}\ln \left(1-\frac{1}{2r}\ln (r)\right)\leq 
\ln\left( \sqrt{\frac{2e}{2e-1}}\right).
\end{equation}
 The estimates \eqref{minmaxvc} of $v_c^{-1}$  follow from \eqref{vcequat}
and \eqref{eee}. Now from   \eqref{minmaxvc}, 
we easily deduce the following asymptotic for $v^{-1}$,
$$
v^{-1}(\rho) \sim \ln\left(\frac{c}{\rho} \sqrt{ \ln(\frac{c}{\rho})}\right)
\sim
\ln \frac{1}{\rho},
$$
when $\rho$ tends to $0$.
Statement (2) is proved.

Note  that the pointwise estimate \eqref{minmaxvc} is  strictly stronger information on 
$v^{-1}(\rho)$ than the equivalent term $\ln (1/\rho)$ due to the second order term $\frac{1}{2}\ln\ln(1/\rho)$ when $\rho$ goes to $0$.
\vskip0.3cm

(3) From Statement (1), for all $\varepsilon \in (0,1)$, there exists $R_{\varepsilon}>1$ such that, 
for all $r\geq R_{\varepsilon}$, we have
$$
g_1(r)\geq v_{c_-}(r),
$$
where $c_{-}=(1- \varepsilon)\sqrt{\frac{\pi}{2}}$.

Let  
$\rho_{\varepsilon}= v_{c_{-}} (R_{\varepsilon})
\leq  \min (g_1 (R_{\varepsilon}), v_{c_{-}} (1)) \leq \min (1, v_{c_{-}} (1))$.
Then, by \eqref{minmaxvc} and the fact that $g_1$ and $v_{c_-}$ are strictly decreasing on $(1,+\infty)$, we have for all $\rho\in (0,\rho_\varepsilon)$,
$$
g_1^{-1}(\rho)\geq v_{c_{-}} ^{-1}(\rho)
\geq 
\ln\left( \frac{{c_{-}} }{\rho}\sqrt{ \ln ( \frac{{c_{-}} }{\rho})}\right).
$$
This proves the lower bound of \eqref{g1inverse}. 
The proof of the upper bound  of  \eqref{g1inverse} is similar and left to the reader. This concludes the proof of Statement (3).

\vskip0.3cm

(4)
From (1), we deduce that there exists $0<c_1<c_2<+\infty$  and $R>1/2$ such that for all $r>R$,
$$
v_1(r)\leq g_1(r)\leq v_2(r),
$$
with $v_i(r)=c_iv(r)$. Since $v_1, v_2$ and $g_1$ are continuous strictly decreasing, 
and goes to zero as $r\rightarrow +\infty$, 
we deduce  that  for all 
$\rho\in (0,\rho_1)$ with  $\rho_1=\min (v_1(R),1)$,
$$
v^{-1}_1(\rho)\leq g^{-1}_1(\rho)\leq v^{-1}_2(\rho),
$$
By (2), we have $v_i^{-1}(\rho)\sim \ln\left( \frac{1}{\rho}\right)$ as $\rho$ goes to $0$ independently of the constants $c_i$. 
This implies that 
$g^{-1}_1(\rho)\sim \ln\left( \frac{1}{\rho}\right)$ as $\rho$ goes to $0$,
 and  proves Statement (4).
 \vskip0.3cm

(5)
This result is  proved by using (1)-(2) of  Proposition \ref{propgun},  and  by the fact that $v^{-1}$  (also $w$) is decreasing.
Here, we just sketch the proof.

By Statement (1) of Proposition \ref{propgun}, we have 
for any  $\varepsilon\in (0,1)$ that  there exists $R_{\varepsilon}>0$ (large enough), such that
\begin{equation}\label{equiv}
(1-\varepsilon)v(r)\leq g_1(r)\leq (1+\varepsilon)v(r),\quad r>R_{\varepsilon},
\end{equation}
where $v(r)=\sqrt{\frac{\pi r}{2}}e^{-r}$. The inverse $v^{-1} $  of $v$ exists and  is a decreasing function
on some interval $(0,\rho_1)$  with  a sufficiently small $\rho_1$. Therefore, we obtain  
$$
v^{-1} 
\left(
(1+\varepsilon)v(r)
\right)
\leq
v^{-1} \circ g_1(r)
\leq 
v^{-1} 
\left(
 (1-\varepsilon)v(r)
\right),
$$
for $r$ large enough. So, it is enough to show that, for any constant $\lambda>0$, we have
\begin{equation}\label{vv}
v^{-1} 
\left(
\lambda  v(r)
\right)\sim r, 
\end{equation}
as  $r\rightarrow +\infty$.
From Statement (2) of Proposition \ref{propgun},
we have $v^{-1}(\rho)\sim \log(1/\rho)$ as $\rho\rightarrow 0$.
Then the asymptotic \eqref{vv} easily follows, since we have
$$
v^{-1} 
\left(
\lambda v(r)
\right)
\sim 
r-
\frac{1}{2}\ln r+K
\sim r,
$$
when $r$  is large enough. (Here, $K$ denotes  some constant depending on $\lambda$, but independent of $r$.)
\vskip0.3cm

\noindent Let  $\varepsilon\in (0,1)$.
Using again \eqref{equiv} with $r=v^{-1}(\rho)$ for small $\rho$ (i.e.  $r$ is large), we obtain
$$
(1-\varepsilon)\rho \leq g_1\circ v^{-1}(\rho)\leq (1+\varepsilon)\rho,
$$
which is exactly  $g_1\circ v^{-1}(\rho)\sim \rho$ as $\rho\rightarrow 0$.
The proof is similar when considering $w$ instead of  $v^{-1}$. 
On one hand, we use the fact that
$w(\rho)\sim\ln(1/\rho)$ as $\rho$ goes to zero.
On the other hand, 
we need careful treatment for  the asymptotics
of  $v\circ w_c(\rho)\sim \rho$ as  $\rho$ goes to $0$.
Indeed, we have by an explicit computation
$$
v\circ w(\rho)
=
\rho \left[1+\frac{1}{2} 
\frac{\ln(\ln \frac{c}{\rho})
}
{ \ln\frac{c}{\rho}}
\right]^{1/2}
\sim \rho, \quad \rho\rightarrow 0,
$$
where $c=\sqrt{\frac{\pi}{2}}$.
(The function $v^{-1}$ can be called an asymptotic inverse of $g_1$). 
\vskip0.3cm

\noindent To conclude the proof of the proposition, note that 
$h(\rho)=\ln (\frac{1}{\rho})$ is not an asymptotic inverse of $g_1$. 
Indeed, we have 
$h(g_1(r))\sim r$ as $r$ goes to infinity but  $g_1(h(\rho))\nsim \rho $ as $\rho$ goes to $0$ due to the fact that
$v\circ h(\rho)=c\rho \sqrt{\ln(1/\rho)}$. 
This last remark   motivates the  introduction of  the function $v_c$ to study the inverse
$g_1^{-1}$.
In particular, this allows us to obtain precise pointwise bounds of $g_1^{-1}$.
 \hfill $\square$
\vskip0.3cm

In the last part of this section, we prove below the asymptotic \eqref{gamgamcc}  mentioned in  Section \ref{aim}. 
To be more specific, we  provide some lower bounds on the growth of the ratio
$\frac{\gamma^*(\alpha,M)}{\gamma^*_{cc}(\alpha,M)}$ when  $M$ goes to $8\pi^+$, from which   \eqref{gamgamcc} can be deduced.

\begin{pro}\label{progamcc}
Let $\gamma^*_{cc}(\alpha,M)$ and $\gamma^*(\alpha,M)$ defined as in the two inequalities
\eqref{gammacc}  and \eqref{eqgamma}
respectively. 
Let $g_{1}$ be defined as in \eqref{g1def}. Then we have
\begin{enumerate}
\item
Let $M>8\pi$. For any  $\beta\in (0, \frac{1}{2})$, there exists a constant $C_{\beta}>0$ such that
\begin{equation}\label{minquot}
 \frac{\gamma^*(\alpha,M)}{\gamma^*_{cc}(\alpha,M)}
 \geq
 C_{\beta}
 \left[
g_{1}^{-1}(8\pi /M)\right]^{-(2-4\beta)}.
\end{equation}
\item
We have
$$
\lim_{M\rightarrow 8\pi^+}  
\frac{\gamma^*(\alpha,M)}{\gamma^*_{cc}(\alpha,M)}=+\infty.
$$
\end{enumerate}
\end{pro}

{\bf Proof}.
The second statement, i.e. \eqref{gamgamcc}, is a  direct consequence of \eqref{minquot}.
Indeed,  the exponent $-(2-4\beta)$ is negative for any  $\beta<1/2$,  and
$\lim_{M\rightarrow 8\pi^+} g_{1}^{-1}(8\pi /M)=g_{1}^{-1}(1)=0$ by continuity of $g_{1}^{-1}$ at $1^{-}$.
\vskip0.3cm

\noindent We prove  \eqref{minquot} as follows. From the definitions of $\gamma^*(\alpha,M)$
and $\gamma^*_{cc}(\alpha,M)$, the following ratio  (independent of $\alpha$)  has the next form
$$
\frac{\gamma^*(\alpha,M)}{\gamma^*_{cc}(\alpha,M)}
=
2{\mathcal C}^2\left[
\frac{ g_{1}^{-1}(8\pi /M)}
{1-8\pi /M}\right]^2.
$$
Let $r=g_{1}^{-1}(8\pi /M)$. Thus,  we have $g_{1}(r)=8\pi /M$, and $M\rightarrow 8\pi^+$ if and only if 
$r\rightarrow 0^+$.
Then, the  ratio  above can  be written in the following form,
$$
\frac{\gamma^*(\alpha,M)}{\gamma^*_{cc}(\alpha,M)}
=
2{\mathcal C}^2\left[
\frac{1- g_{1}(r)}
{r}\right]^{-2}.
$$ 
We now prove that for  any $\beta\in  (0,\frac{1}{2})$ there exists a constant $K_{\beta}>0$ such that, for all $r>0$,
$$
\frac{1- g_{1}(r)}
{r}
\leq 
K_{\beta} \, r^{1-2\beta}.
$$
From which we easily obtain  inequality  \eqref{minquot}.
Indeed, we have the  two  following  integral  representation formula :
$$
\frac{1-g_{1}(r)}{r}=\int_0^{+\infty}\frac{(1-e^{-\frac{r^2}{4t}})}{r}  e^{-t}\, dt,
$$ 
and 
$$
\frac{(1-e^{-\frac{r^2}{4t}})}{r} 
=
\frac{1}{4 tr}
\int_0^{r^2} e^{-s/4t}\, ds.
$$
On the other hand, it is  well-known that for any $\beta>0$  there exists a constant $L_{\beta}=\beta^{\beta}e^{-\beta}>0$ 
(optimal)  such that
$$
e^{-x}\leq 
\frac{L_{\beta}}{x^{\beta}},\quad x>0.
$$
Take  $x=s/4t$ in this inequality  and integrate the corresponding one  with respect to $s$ from $0$ to $r^2$.
For  any $\beta\in (0, 1)$,  $r,t>0$, we thus  get
$$
\frac{(1-e^{-\frac{r^2}{4t}})}{r} 
\leq
\frac{L_{\beta}}{(1-\beta)} 4^{\beta-1} \, t^{\beta-1} r^{1-2\beta}.
$$ 
We finish the proof by multiplying this expression by $e^{-t}$, and integrate  with respect to $t>0$.
We  ultimately  obtain the asserted inequality, i.e.
$$
0<
\frac{1- g_{1}(r)}{r}
\leq
\frac{L_{\beta}}{(1-\beta)} 4^{\beta-1} 
\Gamma(\beta)
\, r^{1-2\beta}
=:
K_{\beta}\, r^{1-2\beta}.
$$
Here,  $\Gamma$ denotes  the classical $gamma$ function.
This completes the proof of Proposition \ref{progamcc}. 
\hfill $\square$


\end{document}